\documentclass[12pt]{article}
\usepackage{amssymb}
\usepackage{latexsym}
\usepackage{fullpage}
\usepackage{amscd}
\def\R{\mathbb{R}}
\def\C{\mathbb{C}}
\def\so{\mathfrak{so}}
\def\sp{\mathfrak{sp}}
\def\SO{\mathbf{SO}\,}
\def\Sp{\mathbf{Sp}\,}

\def\O{\mathbf{O}\,}

\def\hol{\mathfrak{hol}\,}
\def\Hol{\mathrm{Hol}\,}
\def\Hom{\mathrm{Hom}\,}
\def\id{\mathrm{id}\,}
\def\#{\sharp}
\def\b{\flat}

\def\s{\sigma}
\def\x{\mathrm{expt}}
\def\D{\Delta}
\def\L{\Lambda\,}
\def\M{\Lambda_\mathrm{top}}
\def\S{\mathrm{Sym}\,}
\def\T{\mathcal{R}}
\def\W{\mathrm{Schur}\,}
\def\Y{\mathfrak{Y}}
\def\Z{\mathbf{S}\,}
\def\<#1,#2>{\langle#1,#2\rangle}
\def\ind{\mathrm{ind}\,}
\def\rep{\mathop{\bullet}}
\def\inner{\mathop{\raise.1em\hbox{$\lrcorner$}}}

\def\innez{\mathop{\raise.1em\hbox{$\lrcorner_\circ$}}}
\def\qed{\ensuremath{\quad\Box\quad}}
\newtheorem{Proposition}{Proposition}[section]
\newtheorem{Lemma}[Proposition]{Lemma}
\newtheorem{Theorem}[Proposition]{Theorem}
\newtheorem{Corollary}[Proposition]{Corollary}
\newtheorem{Definition}[Proposition]{Definition}

\newtheorem{Remark}[Proposition]{Remark}

\def\proof{\noindent\textbf{Proof:}\quad}
\def\summary#1{(\textsl{#1})\hfill\break}
\def\nosummary{\hfill\break}
\def\pfill{\par\vskip3mm plus1mm minus1mm}
\begin{document}
\title{Vanishing Theorems for Quaternionic K\"ahler Manifolds}
\author{Uwe Semmelmann \& Gregor Weingart
 \footnote{partially supported by the SFB 256
 \textsl{``Nichtlineare partielle Differentialgleichungen''}}}
\maketitle
\begin{abstract}
 In this article we discuss a peculiar interplay between the representation
 theory of the holonomy group of a Riemannian manifold,  the Weitzenb\"ock
 formula for the Hodge--Laplace operator on forms and the Lichnerowicz formula
 for twisted Dirac operators. For quaternionic K\"ahler manifolds this leads
 to simple proofs of eigenvalue estimates for Dirac and Laplace operators.
 Moreover, it enables us to determine which representations can contribute
 to harmonic forms. As a corollary we prove the vanishing of certain odd
 Betti numbers on compact quaternionic K\"ahler manifolds of negative scalar
 curvature. We simplify the proofs of several related  results in the positive
 case.
\end{abstract}
\vskip0.5cm
\tableofcontents
\section{A Prelude on Weitzenb\"ock Formulas}
 
 Since decades the Weitzenb\"ock formulas for Dirac operators on
 Clifford bundles have inspired intensive and important research.
 Beautiful results can be proved elegantly using the full power
 of the Weitzenb\"ock machinery. 
 The basic example of a Clifford bundle is the bundle of exterior
 forms $\L^\bullet T^*M$ endowed with the scalar product induced by
 the metric on $M$ and Clifford multiplication with tangent vectors
 $$
  \star:\;\; T_pM\times\L^\bullet T^*_pM\;\longrightarrow\;
  \L^\bullet T^*_pM,\qquad (X,\,\omega)\;\longmapsto\;X\star\omega
 $$
 defined by $X\star\omega:=X^\#\wedge\omega\,-\,X\inner\omega$. The
 Levi--Civita--connection induces a connection $\nabla$ on $\L^\bullet
 T^*M$ and an associated second order elliptic differential operator
 $\nabla^*\nabla\,:=\,-\sum_i\nabla^2_{E_i,E_i}$ where $\nabla^2_{X,Y}
 \,:=\,\nabla_X\nabla_Y-\nabla_{\nabla_XY}$ and the sum is over a local
 orthonormal base $\{E_i\}$. On the other hand we have the exterior
 differential $d$ and its formal adjoint $d^*$ as natural first order
 differential operators on $\L^\bullet T^*M$ linked to $\nabla^*\nabla$
 by the classical Weitzenb\"ock formula
 \begin{equation}\label{wclass}
  \D\quad:=\quad(d+d^*)^2\quad=\quad\nabla^*\nabla\;+\;{1\over2}
  \sum_{ij}\,E_i\star\,E_j\star\,R_{E_i,E_j}
 \end{equation}
 where $R_{X,Y}$ is the curvature endomorphism of $\L^\bullet T^*_pM$.
 However the connection on $\L^\bullet T^*M$ is induced by a connection
 on $TM$ and consequently the curvature endomorphism $R_{X,Y}$ is just
 the curvature endomorphism of $T_pM$ in a different representation,
 namely the representation
 $$
  \rep:\;\;\so(T_pM)\times\L^\bullet T^*_pM\;\longrightarrow\;
  \L^\bullet T^*_pM,\qquad (X,\,\omega)\;\longmapsto\;X\rep\omega
 $$
 of the Lie algebra $\so(T_pM)$ of $\SO(T_pM)$ on the exterior
 algebra induced by its representation on $T_pM$. The canonical
 identification of $\so(T_pM)$ with the bivectors characterized by
 $$
  \L^2T_pM\;\stackrel\cong\longrightarrow\;\so(T_pM),\qquad
  \<\,(X\wedge Y)\rep A,\,B\,>\;:=\; \<\,X\wedge Y,\,A\wedge B\,>
 $$
 reads $(X\wedge Y)\rep A\;:=\;\<X,A>\,Y\,-\,\<Y,A>\,X$
 and defines a unique bivector $R(X\wedge Y)$ via:
 $$
  \<R(X\wedge Y)\rep Z,W> \quad:=\quad \<R_{X,Y}Z,W> \qquad\qquad
  R(X\wedge Y)\quad=\quad{1\over2}\,\sum_i\,E_i\wedge R_{X,Y}E_i
 $$

 \noindent In the spirit of this identification the representation
 of $\so(T_pM)$ on $\L^\bullet T^*_pM$ is given by $(X\wedge Y)\rep
 \;=\;Y^\#\wedge X\inner\,-\,X^\#\wedge Y\inner$. In particular, the
 classical Weitzenb\"ock formula becomes
 \begin{eqnarray*}
  \D
  &=& \nabla^*\nabla\;+\;{1\over2}\,\sum_{ij}\,
   (E_i^\#\wedge E_j^\#\wedge\,-\,E_i\inner E_j^\#\wedge
   \,-\,E_i^\#\wedge E_j\inner\,+\,E_i\inner E_j\inner)\,
   R(E_i\wedge E_j)\rep \\
  &=& \nabla^*\nabla\;+\;{1\over2}\,\sum_{ij}\,
   (E_i\wedge E_j)\rep R(E_i\wedge E_j)\rep
 \end{eqnarray*}
 because both potentially troublesome inhomogeneous terms cancel
 by the first Bianchi identity leaving us with a curvature term
 depending linearly on the curvature tensor:
 $$
  R\quad:=\quad{1\over 4}\,\sum_{ij}\,(E_i\wedge E_j)\,\cdot\,R(E_i\wedge E_j)
  \qquad\in\quad\S^2(\L^2 T_pM) \ .
 $$

 \noindent
 It will be convenient to compose the identification $\L^2T_pM\,
 \stackrel\cong\longrightarrow\,\so(T_pM)$ with the quantization
 map $\,q:\,\S^2\so(T_pM)\,\longrightarrow\,\mathcal{U}\,\so(T_pM),
 \;X^2\,\longmapsto\,X^2,$ into the universal enveloping algebra
 of $\so(T_pM)$ to get an element $q(R)\in\mathcal{U}\,\so(T_pM)$ with:
 \begin{equation}\label{weitzen}
  \Delta\quad=\quad\nabla^*\nabla\;+\;2\,q(R)
 \end{equation}

 What is the advantage of writing the well known classical Weitzenb\"ock
 formula (\ref{wclass}) in this fancy way? Well, the Weitzenb\"ock
 formula (\ref{weitzen}) brings the holonomy group of the underlying
 manifold into play. Recall that the holonomy group $\Hol_pM\subset\O(T_pM)$
 is the closure of the group of all parallel transports along piecewise
 smooth loops in $p\in M$. We will assume throughout that $M$ is connected
 so that the holonomy groups in different points $p$ and $\tilde p$ are
 conjugated by parallel transport $T_pM\,\longrightarrow\,T_{\tilde p}M$.
 Choosing a suitable representative $\Hol\subset\O_n\R$ with $n:=\dim M$
 of their common conjugacy class acting on the abstract vector space
 $\R^n$ we can define the holonomy bundle of $M$:
 $$
  \Hol(M)\quad:=\quad\{\;f:\;\R^n\,\longrightarrow\,T_pM\,\vert
  \;\;p\in M\textrm{\ and\ }f\textrm{\ isometry with\ }
  f(\Hol)\,=\,\Hol_pM\;\} \ .
 $$
 The holonomy bundle is a reduction of the orthonormal frame bundle
 $\O(M)$ to a principal bundle with structure group $\Hol$, which is
 stable under parallel transport. Consequently the Levi--Civita
 connection is tangent to $\Hol(M)$ and descends to a connection
 on $\Hol(M)$.

 The associated fibre bundle $\Hol(M)\times_\Hol\O_n\R$ is
 canonically diffeomorphic to the full orthonormal frame bundle
 $\O(M)$. This construction provides an explicit foliation of
 $\O(M)$ into mutually equivalent principal subbundles stable
 under parallel transport. Choosing a leaf different from the
 distinguished leaf $\Hol(M)$ amounts to choosing a different
 representative for the conjugacy class of $\Hol\subset\O_n\R$.
 In particular every principal subbundle of $\O(M)$ stable under
 parallel transport is a union of leaves and is characterized
 by a subgroup of $\O_n\R$ containing a representative of the
 conjugacy class of the holonomy group $\Hol$.

 With the Levi--Civita connection being tangent to the holonomy
 bundle $\Hol(M)$ its curvature tensor $R$ takes values in the
 holonomy algebra $\hol_pM$ at every point $p\in M$, so that
 $R\,\in\,\S^2\hol_pM\,\subset\,\S^2\L^2T_pM$ and $q(R)\,\in\,
 \mathcal{U}\,\hol_pM$. However by definition every point
 $f\,\in\,\Hol(M)$ identifies $\hol_pM$ with $\hol$ making
 $q(R)$ a $\mathcal{U}\,\hol$--valued function on $\Hol(M)$:
 $$
  q(R)\quad\in\quad C^\infty(\;\Hol(M),\,\mathcal{U}\,\hol\;)^{\Hol}
   \;\cong\;\Gamma(\Hol(M)\times_\Hol\mathcal{U}\,\hol) 
 $$
 For an arbitrary irreducible complex representation $\pi$ of
 $\Hol$ the associated vector bundle $\pi(M)\,:=\,\Hol(M)\times_{\Hol}\pi$
 over $M$ is endowed with the connection induced from the Levi--Civita
 connection. Moreover there is a canonical second order differential
 operator defined on sections of $\pi(M)$:
 \begin{equation}\label{dpi}
  \D_\pi\quad:=\quad \nabla^*\nabla\;+\;2\,q(R)
 \end{equation}
 
 \noindent It is evident from the Weitzenb\"ock formula (\ref{wclass})
 written as in (\ref{weitzen}) that the diagram
 $$
  \begin{CD}
   \pi(M) @>{\textstyle\D_\pi}>> \pi(M) \\
   @V{\textstyle F}VV @VV{\textstyle F}V \\
   \L^\bullet T^*M\otimes_\R\C @>{\textstyle\D}>>\L^\bullet T^*M\otimes_\R\C
  \end{CD}
 $$
 commutes for any $F\,\in\,\Hom_{\Hol}(\pi,\L^\bullet\C^{n*})$ or
 equivalently for any globally parallel embedding $F:\,\pi(M)\,
 \longrightarrow\,\L^\bullet T^*M\otimes_\R\C$. Hence the pointwise
 decomposition of $\L^\bullet T^*_pM\otimes_\R\C$ into irreducible
 complex representations of $\Hol_pM$ becomes a global decomposition
 of any eigenspace of $\D$, e.~g.~we have for its kernel:
 $$
  H_{dR}^\bullet(M,\,\C)\;\;=\;\;\bigoplus_\pi\Hom_{\Hol}(\,\pi,
  \,\L^\bullet \C^{n*}\,)\,\otimes\,\mathrm{Kern}\;\D_\pi
 $$

 The same kind of reasoning is possible for the Dirac operator on
 spinors, assuming the manifold $M$ to be spin and taking $\Hol_pM$
 to be its spin holonomy group. Ignoring for the moment the Lichnerowicz
 result that the curvature term reduces to multiplication by the scalar
 curvature and employing the formula $(X\wedge Y)\rep\,:=\,{1\over2}
 (\,X\star Y\star\,+\,\<X,Y>\,)$ for the representation of $\so(T_pM)$ on
 the spinor bundle $\Z(M)$ we can proceed from (\ref{wclass}) directly to:
 \begin{equation}\label{dirac}
  D^2\quad=\quad\nabla^*\nabla\;+\;4\,q(R)\,.
 \end{equation}

 \noindent
 In particular, all eigenspaces of $D^2$ decompose globally according to
 the pointwise decomposition of the spinor bundle under the spin holonomy
 group $\Hol_pM$. Whereas the change of the factor of $q(R)$ from 2 to 4
 is certainly puzzling, there can be no doubt however that equation
 (\ref{dirac}) is true. In fact from Lichnerowicz's result we already
 know that $q(R)$ acts by scalar multiplication with ${\kappa\over 16}$
 on $\Z(M)$, where $\kappa$ is the scalar curvature of $M$. Hence we can
 read equation (\ref{dirac}) as
 $$
  D^2\Big\vert_\pi \quad=\quad \Delta_\pi\;+\;{\kappa\over8}
 $$
 where the restriction to $\pi$ is a short hand notation for any 
 globally parallel embedding $F:\;\pi(M)\,\longrightarrow\,\Z(M)$ induced
 by some non--trivial $F\,\in\,\Hom_{\Hol}(\pi,\Z)$. Written in this way
 formula (\ref{dirac}) is seen to be a generalization of the Partharasathy
 formula for the Dirac square $D^2$ on a symmetric space $G/K$ of compact
 type, because in this case the operators $\Delta_\pi$ defined above on
 sections of $\pi(M)$ all become the Casimir of $G$.

 At this point the reader may argue that these results are not too
 surprising because intrinsically defined differential operators are
 restricted to parallel subbundles. However the main point is that
 $\D$ and $D^2$ do not only respect some decomposition into parallel
 subbundles, but that their restrictions to these subbundles are
 completely independent of the embedding. Counterexamples to the
 idea that intrinsically defined differential operators always
 enjoy these two properties are easily found among twisted Dirac
 operators.

 Consider therefore a geometric vector bundle $\T(M)\,:=\,\Hol(M)
 \times_\Hol\T$ associated to the holonomy bundle via some not
 necessarily irreducible representation $\T$ of the holonomy group.
 The Levi--Civita connection on $\Hol(M)$ defines a geometric connection
 on this vector bundle, whose curvature endomorphism is still given
 through the representation
 $$
  \rep:\;\;\hol_pM\times\T_p(M)\;\longrightarrow\;\T_p(M)
 $$
 of the Lie algebra $\hol_pM$ on $\T_p(M)$ by the formula
 $R_{X,Y}^\T\,=\,R(X\wedge Y)\rep$. The twisted Dirac operator
 $D_\T$ is a first order differential operator acting on sections
 of the vector bundle $(\Z\otimes\T)(M)$. It satisfies a twisted
 Weitzenb\"ock formula derived from (\ref{wclass}):
 \begin{equation}\label{twist}
  D^2_\T\quad=\quad\nabla^*\nabla\;+\;{1\over2}
  \sum_{ij}\,\Big(E_i\star\,E_j\star\,R(E_i\wedge E_j)\rep\,\otimes\,
  \id_\T\;+\;E_i\star\,E_j\star\,\otimes R(E_i\wedge E_j)\rep\Big)
 \end{equation}
 
 \noindent
 This formula has an apparent asymmetry between the spinor bundle and
 the twist. However, we still have the formula $(X\wedge Y)\rep\,=\,
 {1\over2}(\,X\star\,Y\star\,+\,\<X,Y>\,)$ for the representation of
 $\so(T_pM)$ on the fibre $\Z_p(M)$  of the spinor bundle and we may try to 
 balance this asymmetry to cast equation (\ref{twist}) into a form similar to
 (\ref{dirac}). This is most easily achieved by rewriting the action of
 $q(R)$ on the tensor product $\Z\otimes\T$ in the following asymmetric way:
 \begin{eqnarray}
  q(R)&=& {1\over2}\sum_{ij}\Big(\;
    (E_i\wedge E_j)\rep R(E_i\wedge E_j)\rep\otimes\id_\T\;+\;
    (E_i\wedge E_j)\rep \otimes R(E_i\wedge E_j)\rep\;\Big)\\[-2pt]
  && \qquad\qquad -\;\;q(R)\otimes\id_\T\;\;+\;\;
     \id_\Z\otimes\,q(R)\nonumber
 \end{eqnarray}
 With Lichnerowicz's result $q(R)\,=\,{\kappa\over16}$ for
 the spinor representation $\Z$ equation (\ref{twist}) becomes
 \begin{equation}\label{tdir}
  D^2_\T\quad=\quad\D_{\Z\otimes\T}
   \;+\;{\kappa\over8}\,\otimes\,\id_\T
   \;-\;\id_\Z\,\otimes\,2\,q(R)
 \end{equation}

 \noindent
 In conclusion, the squares $D^2_\T$ of twisted Dirac operators will in
 general not respect the decomposition of $(\Z\otimes\T)(M)$ into parallel  
 subbundles because of the critical summand $\;\id_\Z\,\otimes\,2\,q(R)$.
 Nevertheless, if $q(R)$ acts by scalar multiplication not only on $\Z$ but on
 $\T$, too, the global decomposition of the eigenspaces of $D^2_\T$
 according to the pointwise decomposition of $\Z\otimes\T$ is restored. 

 Equation (\ref{tdir}) is the key relation of this article and forms the
 cornerstone and motivation of all statements and calculations to come.
 In fact, we can take advantage of equation (\ref{tdir}) even if the
 manifold in question is not spin, because the twisted Dirac operator may
 be well defined on the vector bundle $(\Z\otimes\T)(M)$ although $M$ is
 neither spin nor $\Z(M)$ or $\T(M)$ are well defined vector bundles.
 The only thing that really matters is whether the representation
 $\Z\otimes\T$ is defined for the holonomy group $\Hol$ itself or
 only for some covering group.

\section{Quaternionic K\"ahler Holonomy}

 In this section we introduce the main notions of quaternionic K\"ahler
 holonomy based on the group $\Hol\,=\,\Sp(1)\cdot\Sp(n)$ with $n\,\geq\,2$.
 Very few examples of compact manifolds with this particular holonomy
 group are known, and it is a deep result that in every quaternionic 
 dimension $n$ there are up to isometry only finitely many of these 
 manifolds with positive scalar curvature $\kappa\,>\,0$ (\cite{lebrun}).
 In fact, the only known examples with $\kappa\,>\,0$ are symmetric spaces,
 the so-called Wolf spaces.

 In order to introduce quaternionic K\"ahler holonomy we return for a
 moment to a point we glossed over in the definition of the holonomy bundle.
 There we had to choose a suitable representative $\Hol\,\subset\,\O_{4n}\R$
 in the conjugacy class of the holonomy groups acting on an abstract vector
 space $\R^{4n}$. This abstract vector space has no meaning in itself but
 plays the role of the tangent representation of $\Hol$ just as $T_pM$ is
 the tangent representation of $\Hol_pM$. Instead of really choosing a
 representative $\Hol\,\subset\,\O_{4n}\R$ it is always better to start
 with specifying this tangent representation. Let us begin with an abstract
 complex vector space $E\,\cong\,\C^{2n}$ endowed with a symplectic form
 $\s\,\in\,\L^2E^*$ and an adapted, positive quaternionic structure $J$,
 i.~e.,~a conjugate linear map $J:\;E\,\longrightarrow\,E$ satisfying
 $$
  J^2\quad=\quad-1 \qquad\qquad
  \s(Je_1,\,Je_2)\quad=\quad\overline{\s(e_1,\,e_2)}\qquad\qquad
  \s(e,\,Je)\quad>\quad 0
 $$
 for all $e_1,\,e_2\,\in\,E$ and $e\,\neq\,0$. Such a set of structures
 is consistent and can be defined on the underlying complex vector space
 of $\mathbb{H}^n$. One merit of this explicit construction is that the
 group of all symplectic transformations of $E$ commuting  with $J$ agrees
 in this picture with the quaternionic unitary group $\Sp(n)\,:=\,\{\,A\in
 \mathrm{M}_{n,\,n}\mathbb{H}\textrm{\ such that\ }\overline{A}^tA=1\,\}$.
 The symplectic form $\s$ induces mutually inverse isomorphisms $\#:\;E\,
 \longrightarrow\,E^*,\,e\,\longmapsto\,\s(e,\cdot)$ and $\b:\;E^*\,
 \longrightarrow\,E$. Similar to the representation of $\L^2T_pM$
 on $T_pM$ considered in the first section there is an action
 $$
  \rep:\;\;\S^2\,E\times E\;\longrightarrow\;E,\qquad (e_1e_2,\,e)\;
  \longmapsto\;(e_1e_2)\rep e\,:=\,\s(e_1,\,e)e_2\,+\,\s(e_2,\,e)e_1
 $$
 of the second symmetric power $\S^2E$ on $E$. This action is skew symplectic
 and commutes with $J$ for all real elements of $\S^2E$. It identifies this
 real subspace with the Lie algebra $\sp(n)$ of $\Sp(n)$ and makes $\rep$
 not only an action but a representation.

 Let $H\,\cong\,\C^2$ be another abstract vector space with the same
 structures, a symplectic form $\s\,\in\,\L^2H^*$ and an adapted, positive
 quaternionic structure $J$. The tensor product $H\otimes E$ of these two
 vector spaces carries a real structure $J\otimes J$ and a complex bilinear
 symmetric form $\<\,,\,>\,:=\,\s\otimes\s$, which is positive definite on
 the real subspace. In this way the group $\O(H\otimes E)$ of all complex
 linear isometries of $H\otimes E$ commuting with $J\otimes J$ is isomorphic
 to $\O_{4n}\R$ and has a distinguished subgroup $\Sp(1)\cdot\Sp(n)\;:=\;
 \Sp(1)\times\Sp(n)/\mathbb{Z}_2$ preserving the tensor product structure
 of $H\otimes E$:

 \begin{Definition}\label{qkh}
 \summary{Quaternionic K\"ahler Manifolds}
  A quaternionic K\"ahler manifold $M$ is a Riemannian manifold of
  dimension $4n,\,n\geq 2,$ endowed with a reduction of the frame
  bundle $\O(M)$ to a principal $\Sp(1)\cdot\Sp(n)$--bundle
  $\Sp(1)\cdot\Sp(M)$ stable under parallel transport. Such a
  reduction exists if and only if the holonomy group $\Hol$
  of $M$ is conjugated to a subgroup of $\Sp(1)\cdot\Sp(n)\,
  \subset\,\O(H\otimes E)$ and in case of equality it may
  be defined as:
  $$
   \Sp(1)\cdot\Sp(M):=
   \{\,f:H\otimes E\longrightarrow T_pM\otimes_\R\C\,\vert\,
    f\textrm{\ isometry and\ }f(\Sp(1)\cdot\Sp(n))=\Hol_pM\,\}\,.
  $$
  If the holonomy group of a quaternionic K\"ahler manifold
  $M$ is conjugated to a proper subgroup of $ \Sp(1)\cdot\Sp(n)$,
  then $M$ is necessarily locally symmetric and its universal
  covering is a Wolf space.
 \end{Definition}

 There are a few remarks to make on this definition. First of all
 we insist on $n\,\geq\,2$, because taking this definition as it stands
 it applies to every oriented Riemannian manifold $M$ of dimension 4.
 In addition a quaternionic K\"ahler manifold with vanishing scalar
 curvature $\kappa\,=\,0$ is locally hyperk\"ahler, its universal
 cover thus hyperk\"ahler, and we will usually exclude these manifolds
 from consideration. In general, however, a quaternionic K\"ahler manifold
 with non--vanishing scalar curvature is despite nomenclature not K\"ahler.

 In order to justify terminology after all these negative remarks and to
 get into contact with a more common definition of quaternionic K\"ahler
 manifolds we recall that $\S^2H$ acts via $(h_1h_2)\rep h\,:=\,
 \s(h_1,\,h)h_2\,+\,\s(h_2,\,h)h_1$ on $H$. For a normed real element
 $i\,h\,Jh\in\S^2H$ with $\s(h,\,Jh)=1$ the action on $H$
 commutes with $J$ and satisfies:
 $$
  (i\,h\,Jh)\rep(i\,h\,Jh)\rep
  \quad=\quad -\id \ .
 $$
 This follows from the fundamental identity $\s(h_1,\,h)h_2\,-\,\s(h_2,\,h)
 h_1\,=\,\s(h_1,\,h_2)h$ for $2$--dimensional symplectic vector spaces and
 hence does not work for $E$. Extending this action from $H$ to the tangent
 representation $H\otimes E$ we conclude that normed real local sections
 of the parallel subbundle $\Sp(1)\cdot\Sp(M)\times_{\Sp(1)\cdot\Sp(n)}\S^2H$
 of the complexified endomorphism bundle $\mathrm{End}\,(TM\otimes_\R\C)$
 act as local complex structures on the tangent bundle $TM$. Choosing
 in this way three local complex structures $I,\,J$ and $K$ satisfying
 $IJ\,=\,K$ we define the canonical quaternionic orientation of $M$
 by declaring every base of the form $X_1,\,IX_1,\,JX_1,\,KX_1,\,\ldots,
 \,X_n,\,IX_n,\,JX_n,\,KX_n$ to be positively oriented. Alternatively the
 canonical quaternionic orientation is induced by the $n$--th power of the
 parallel Kraines form $\Omega\in\L^4(H\otimes E)$ defined in (\cite{krain}).

 A rather subtle remark concerns the two representations $H$ and $E$,
 which do not factor through the projection $\Sp(1)\times\Sp(n)\,
 \longrightarrow\,\Sp(1)\cdot\Sp(n)$. Although we may think of the
 complex tangent bundle as a tensor product of two complex vector
 bundles $H$ and $E$, these vector bundles are not well defined and
 in general exist only locally. In passing from representation theory
 to geometry we always have to check, whether the representations
 factor through the projection $\Sp(1)\times\Sp(n)\,\longrightarrow\,
 \Sp(1)\cdot\Sp(n)$. Things get actually simpler in some respect, as
 the spinor representation $\Z$ of $\Sp(1)\times\Sp(n)$ factors through
 to a representation of $\Sp(1)\cdot\Sp(n)$ whenever $n$ is even. Thus
 all quaternionic K\"ahler manifolds of even quaternionic dimension $n$
 are spin:

 \begin{Proposition}\label{dez}
 \summary{The Signed Spinor Representation (\cite{barsal},\ \cite{wang})}
  The spinor representation $\Z$ of $\Sp(1)\times\Sp(n)$ decomposes
  into the direct sum
  \begin{equation}\label{deco}
   \Z\quad=\quad\bigoplus^n_{r=0}\;\Z_r
     \quad:=\quad\bigoplus^n_{r=0}\;\S^r H\otimes\L^{n-r}_\circ E
  \end{equation}
  where $\L^{n-r}_\circ E$ is the kernel of the contraction
  $\s:\;\L^{n-r}E\,\longrightarrow\,\L^{n-r-2}E$ with the symplectic
  form. For the canonical quaternionic orientation of $H\otimes E$
  the half spin representations are given by:
  $$
   \Z^+\quad:=\quad\bigoplus_{r\equiv n\,(2)}\;\Z_r
   \qquad\qquad
   \Z^-\quad:=\quad\bigoplus_{r\not\equiv n\,(2)}\;\Z_r\,.
  $$
 \end{Proposition}

 The delicate point in an explicit proof of this proposition avoiding
 representation theory is the choice of Clifford multiplication $\star:\,
 (H\otimes E)\times\Z\,\longrightarrow\,\Z$. Besides the Clifford identity
 \begin{equation}\label{clifid}
  (h_1\otimes e_1)\star\,(h_2\otimes e_2)\star\,\,+
  \,(h_2\otimes e_2)\star\,(h_1\otimes e_1)\star\,\quad=\quad
  -\;2\,\s(h_1,\,h_2)\,\s(e_1,\,e_2)
 \end{equation}
 which has to be satisfied, there is another crucial property of this
 multiplication, namely the compatibility condition with the action of
 the Lie algebra $\sp(1)\oplus\sp(n)$ on $\Z$. The representation
 $\bullet$ of the complexified Lie algebra $\S^2H\oplus\S^2E$ of
 the group $\Sp(1)\times\Sp(n)$ on $\Z$ has to agree with the
 representation implicitly defined by Clifford multiplication
 via $(X\wedge Y)\rep\,:=\,{1\over2}(X\star\,Y\star\,+\,\<X,Y>)$.
 This condition depends on the correct formulation of the embedding
 $\S^2H\oplus\S^2E\,\longrightarrow\,\L^2(TM\otimes_\R\C)$. Choosing
 dual pairs of bases $\{de_\mu\},\,\{e_\nu\}$ for $E^*,\,E$ with
 $\<de_\mu,\,e_\nu>\,=\,\delta_{\mu\nu}$ and $\{dh_\alpha\},\,
 \{h_\beta\}$ for $H^*,\,H$ we can check that

 \begin{equation}\label{so}
  (e\,\tilde e)\;\;\longmapsto\;\;
  \sum_\alpha\;(dh_\alpha^\b\otimes e)\wedge(h_\alpha\otimes \tilde e)
  \qquad\qquad
  (h\,\tilde h)\;\;\longmapsto\;\;
  \sum_\mu\; (h\otimes de_\mu^\b)\wedge(\tilde h\otimes e_\mu)
 \end{equation}

 \noindent
 is the correct choice intertwining the representations of $\S^2H,\,\S^2E$
 and $\L^2(TM\otimes_\R\C)$ on $H\otimes E\,=\,TM\otimes_\R\C$. Consequently
 the following two operator identities on the spinor representation $\Z$
 are at the heart of Proposition \ref{dez}:
 \begin{eqnarray}
  (e\,\tilde e)\rep&=&{1\over2}\;\sum_\alpha\,\Big((dh_\alpha^\b\otimes e)
   \star(h_\alpha\otimes\tilde e)\star\,+\,\s(e,\,\tilde e)\Big)\label{ei}\\
  (h\,\tilde h)\rep&=&{1\over2}\;\sum_\mu\,\Big((h\otimes de_\mu^\b)\star
   (\tilde h\otimes e_\mu)\star\,+\,\s(h,\,\tilde h)\Big)\label{hi}
 \end{eqnarray}

 \noindent
 We will not go into the details of this construction given in \cite{qk1},
 but will take Proposition~\ref{dez} as the assertion that a Clifford
 multiplication $\star:\;(H\otimes E)\times\Z\,\longrightarrow\,\Z$
 with the properties (\ref{ei}) and (\ref{hi}) exists satisfying the
 Clifford identity (\ref{clifid}). 

 The most important point in our present discussion of quaternionic
 K\"ahler holonomy is of course the discussion of the curvature tensor
 of a quaternionic K\"ahler manifold and of the associated element
 $q(R)$ in the universal enveloping algebra of the Lie algebra
 $\sp(1)\oplus\sp(n)$ of the holonomy group $\Sp(1)\cdot\Sp(n)$. In fact
 compared to other holonomy groups quaternionic K\"ahler holonomy is
 rather rigid. This is mainly due to the fact that the curvature tensor
 of a quaternionic K\"ahler manifold has to satisfy very stringent
 constraints and can be described completely by the scalar curvature
 $\kappa$ and a section $\mathfrak{R}$ of $\S^4E$. This decomposition
 was first derived by D.~V.~Alekseevskii (cf.:\cite{alex3} or\cite{sala}) 
 and can be made explicit in the following way:

 \begin{Lemma}\label{curvature}
 \summary{The Curvature Tensor}
  A quaternionic K\"ahler manifold $M$ is Einstein with constant
  scalar curvature $\kappa$. Its curvature tensor depends only
  on $\kappa$ and a section $\mathfrak{R}$ of $\S^4E$, this dependence
  reads
  \begin{equation}
   R \quad=\quad -{\kappa\over 8n(n+2)}(R^H + R^E) + R^{hyper}
  \end{equation}
  where the endomorphism valued two forms $R^H,\,R^E$
  and $R^{hyper}$ are defined by:
  \begin{eqnarray}
   R^H_{h_1\otimes e_1,h_2\otimes e_2} &=&
    \s_E(e_1,e_2)(h_1h_2\rep \otimes \id_E) \nonumber\\
   R^E_{h_1\otimes e_1,h_2\otimes e_2} &=&
    \s_H(h_1,h_2)(\id_H \otimes e_1e_2\rep)\\
   R^{hyper}_{h_1\otimes e_1,h_2\otimes e_2} &=& \s_H(h_1,h_2)
   (\id_H \otimes(e_2^\#\inner\,e_1^\#\inner\mathfrak{R})\rep)\nonumber
  \end{eqnarray}
 \end{Lemma}
 
 We will give a short sketch of the proof of this lemma, but refrain from 
 giving all the details which again can be found in \cite{qk1}. The essential 
 point is to show that the linear space of $\Sp(1)\cdot\Sp(n)$--curvature 
 tensors, i.~e.,~the intersection of $\S^2\hol\,\subset\,\S^2\L^2TM$ with
 the kernel of the Bianchi identity $\wedge:\;\S^2\L^2TM\,\longrightarrow
 \,\L^4TM$, is isomorphic to $\C\oplus\S^4E$. Consequently our ansatz for
 $R$ as a linear combination of $R^H\,+\,R^E$ and $R^{hyper}$ is justified
 since $R^H\,+\,R^E$ and $R^{hyper}$ separately satisfy the first Bianchi
 identity. Note that $(e_2^\#\inner\,e_1^\#\inner\,\mathfrak{R})\rep e\,=
 \,\mathfrak{R}(e_1^\#,\,e_2^\#,\,e^\#,\,\cdot)$ is symmetric in $e_1,\,e_2$
 and $e$. In order to determine the curvature tensor $R$ completely, it is
 convenient to calculate the Ricci curvature of $M$, given by the trace
 $\mathrm{Ric}(X,Y) = \mathrm{tr}(Z \longmapsto R_{Z,X}Y)$ of the
 endomorphism $Z\longmapsto R_{Z,X}Y$. The different tensors
 contribute to this endomorphism in the following way:
 \begin{eqnarray*}
  R^H_{h\otimes e,\,h_1\otimes e_1}\,h_2\otimes e_2 &=&
   (\,\s(h,\,h_2)h_1\,+\,\s(h_1,\,h_2)h\,)\otimes \s(e,\,e_1)e_2\\[1pt]
  R^E_{h\otimes e,\,h_1\otimes e_1}\,h_2\otimes e_2 &=&
   \s(h,\,h_1) h_2\otimes(\,\s(e,\,e_2)e_1\,+\,\s(e_1,\,e_2)e\,)\\[1pt]
  R^{hyper}_{h\otimes e,\,h_1\otimes e_1}\,h_2\otimes e_2 &=&
   \s(h,\,h_1) h_2\otimes \mathfrak{R}(e^\#,\,e_1^\#,\,e_2^\#,\,\cdot)
 \end{eqnarray*}
 
 \noindent
 Note that all these endomorphisms preserve the tensor product structure.
 Hence their traces are the product of the partial traces in each tensor
 factor. However, $e\,\longmapsto\mathcal{R}(e_1^\#,\,e_2^\#,\,e^\#,\,\cdot)$
 is induced by an element of $\S^2E$ and hence trace--free, which rules out
 contributions from $R^{hyper}$ to the Ricci curvature. As the trace of the
 endomorphism $e\,\longmapsto\,\s(e,\,e_2)e_1$ is $\s(e_1,\,e_2)$ the
 trace of $e\,\longmapsto\,\s(e,\,e_2)e_1\,+\,\s(e_1,\,e_2)e$ is given
 by $(2n+1)\,\s(e_1,\,e_2)$. Similar remarks apply to $H$ and we are left
 with:

 \begin{equation}\label{ricci}
  (\mathrm{Ric}^H\,+\mathrm{Ric}^E)(\,h_1\otimes e_1,\,h_2\otimes e_2\,)
  \quad=\quad
  -\,(\,2n\,+\,4\,)\,\s(h_1,\,h_2)\,\s(e_1,\,e_2)
 \end{equation}
 
 \noindent
 The Ricci curvature being a multiple of the metric the quaternionic
 K\"ahler manifold $M$ is Einstein, a fortiori the scalar curvature
 $\kappa$ is constant on $M$ and equation (\ref{ricci}) fixes the
 coefficient of $R^H+R^E$ in $R$ via $\mathrm{Ric}(X,\,Y)
 \,=\,{\kappa\over 4n}\<X,\,Y>$.
 \pfill

 At the end of this section we want to describe the action of the
 element $q(R)$ of the universal enveloping algebra $\mathcal{U}
 (\,\sp(1)\oplus\sp(n)\,)$ on some representations. In particular
 we will see that for a large class of representations of $\Sp(1)
 \times\Sp(n)$ the element $q(R)$ acts by scalar multiplication,
 because the contributions from the hyperk\"ahler part $R^{hyper}$
 of the curvature tensor drop out. Observe first that $q(R)$ depends
 linearly on $R$:
 $$
  q(R)\quad=\quad -{\kappa\over 8n(n+2)}\,\left(\,q(R^H)
  \,+\,q(R^E\,)\right)\;+\;q(R^{hyper})
 $$
 Using equation (\ref{so}) we can write down the terms
 appearing in this sum  more explicitly:

 \begin{Lemma}\label{three}
 \begin{eqnarray*}
  q(R^H)\;\;&=&{1\over4}\,\sum_{\alpha\beta}\,
  (dh^\b_\alpha\,dh^\b_\beta)\rep\,(h_\alpha\,h_\beta)\rep\\
  q(R^E)\;\;&=&{1\over4}\,\sum_{\mu\nu}\,
  (de^\b_\mu\,de^\b_\nu)\rep\,(e_\mu\,e_\nu)\rep\\
  q(R^{hyper})&=&{1\over4}\,\sum_{\mu\nu}\,
  (de^\b_\mu\,de^\b_\nu)\rep\,(e_\mu^\#\inner\,e_\nu^\#\inner\,
  \mathfrak{R})\rep
 \end{eqnarray*}
 \end{Lemma}

 \proof
 Converting the sum over a local orthonormal base $\{E_i\}$ into the sum
 $$
  \sum_i\;E_i\otimes E_i\quad=\quad
  \sum_{\alpha\mu}\;(dh_\alpha^\b\otimes de_\mu^\b)
   \otimes (h_\alpha\otimes e_\mu)
 $$
 over dual pairs $\{de_\mu\},\,\{e_\mu\}$ and $\{dh_\alpha\},\,
 \{h_\alpha\}$ of bases we calculate say for $q(R^{hyper})$
 \begin{eqnarray*}
  {1\over4}\,\sum_{ij}\,(E_i\wedge E_j)\rep\,R^{hyper}_{E_i,\,E_j}
  &=&
  {1\over4}\,\sum_{\alpha\beta\mu\nu}\,
   (dh^\b_\alpha\otimes de^\b_\mu\,\wedge\,dh^\b_\beta\otimes de^\b_\nu)\rep
   \,\s(h_\alpha,\,h_\beta)(e_\mu^\#\inner\,e_\nu^\#\inner\,
   \mathfrak{R})\rep\\
  &=&
  {1\over4}\,\sum_{\alpha\mu\nu}\,(dh^\b_\alpha\otimes de^\b_\mu\,
  \wedge\,h_\alpha\otimes de^\b_\nu)\rep\,(e^\#_\mu\inner\,
  e^\#_\nu\inner\,\mathfrak{R})\rep
 \end{eqnarray*}
 which is equivalent to the stated equality in view of equation (\ref{so}).
 \qed
 \pfill

 Evidently $2q(R^H)$ and $2q(R^E)$ respectively are the Casimir operators
 for $\sp(1)$ and $\sp(n)$ in $\s$--normalization, i.~e.,~when the defining
 invariant symmetric form on the Lie algebra $\S^2H$ or $\S^2E$ is not the
 Killing form itself but the natural extension of $\s$ to the second symmetric
 powers using Gram's permanent. For some simple irreducible representations
 it is easy to calculate the Casimir eigenvalues of $q(R^H)$ and $q(R^E)$
 directly. Strictly speaking this procedure is unnecessary because the
 general formula for these eigenvalues in terms of the highest weight
 is simple enough. In this way, however, we get all the Casimir eigenvalues
 we will need below and the precise relations to the Casimirs in Killing
 normalization:

 \begin{Lemma}\label{casim}
 \summary{Casimir Eigenvalues}
  For the irreducible representations $\S^lE$ and $\L^d_\circ E$
  the Casimir eigenvalues for $q(R^E)$ are:
  $$
   q(R^E)_{\S^lE}\;\;=\;\; -\,l\,(\,n\,+\,{l\over2}\,) \qquad\qquad\qquad
   q(R^E)_{\L^d_\circ E}\;\;=\;\; -\,d\,(\,n\,-\,{d\over2}\,+\,1\,)
  $$
 \end{Lemma}

 \proof
 Both calculations are very similar, for the symmetric power $\S^lE$ we get:
 \begin{eqnarray*}
 {1\over4}\sum_{\mu\nu}\,(de^\b_\mu\,de^\b_\nu)\rep\,(e_\mu\,e_\nu)\rep\,
  &=&{1\over4}\sum_{\mu\nu}\,
     (de^\b_\nu\cdot\,de_\mu\inner\,+\,de^\b_\mu\cdot\,de_\nu\inner)\,
     (e_\nu\cdot\,e^\#_\mu\inner\,+\,e_\mu\cdot\,e^\#_\nu\inner)\\
  &=&{1\over2}\sum_{\mu\nu}\,
     (de^\b_\nu\cdot\,\delta_{\mu\nu}\;e^\#_\mu\inner\,+\,de^\b_\nu\cdot\,
     e^\#_\nu\inner\,+\,de^\b_\nu\cdot\,e_\mu\cdot\,de_\mu\inner\,
     e^\#_\nu\inner)\\
  &=&{1\over2}(\,-\,l\,-\,2nl\,-\,l(l\,-\,1)\,)
     \quad=\quad-\,l\,(n\,+\,{l\over2})\qed
 \end{eqnarray*}

 The eigenvalues of $q(R^H)$ are given by the same formulas with
 $n=1$. Setting $l=2$ we get the Casimir eigenvalues for $q(R^E)$
 and $q(R^H)$ in the adjoint representations $\S^2E$ and $\S^2H$ of
 $\sp(n)$ and $\sp(1)$. Since by definition the Casimir eigenvalue
 of the adjoint representation is always one for Casimirs in the
 Killing normalization we get:
 $$
  q(R^E)\quad=\quad-\,2\,(n\,+\,1)\,\mathrm{Cas}_{\sp(n)}
  \qquad\qquad\qquad q(R^H)\quad=\quad-\,4\,\mathrm{Cas}_{\sp(1)}
 $$

 Now we claim that the hyperk\"ahler contribution $q(R^{hyper})$ to the
 element $q(R)$ acts trivially on every irreducible representation occurring
 in the representation $\L E$, i.~e.,~on all representations $\L^d_\circ E$
 with $d\,=\,0,\,\ldots,\,n$. Because $q(R^{hyper})$ depends linearly on
 $\mathfrak{R}\,\in\,\S^4E$ we are allowed to expand $\mathfrak{R}$ into
 a sum of fourth powers ${1\over24}e^4,\,e\in E$, to calculate $q(R^{hyper})$.
 It is thus sufficient to prove that the action of $q({1\over24}e^4)$ on
 $\L E$ is trivial for all $e\,\in\,E$. According to Lemma \ref{three}
 the element $q({1\over24}e^4)$ acts on $\L E$ as:
 $$
  q({1\over24}e^4)
  \;\;=\;\;{1\over2}\,({1\over2}e^2)\rep\,({1\over2}e^2)\rep
  \;\;=\;\;{1\over2}\,(e\wedge\,e^\#\inner)\,(e\wedge\,e^\#\inner)
  \;\;=\;\;-\,{1\over2}\,e\wedge\,e\wedge\,e^\#\inner\,e^\#\inner
  \;\;=\;\;0\,.
 $$
 
 \noindent
 Consequently the curvature tensor $q(R)$ will act by scalar multiplication
 on all representations $\T^{l,\,d}\,:=\,\S^lH\otimes\L^d_\circ E$. From
 equation (\ref{tdir}) we conclude that the squares $D^2_{\T^{l,\,d}}$ of
 the twisted Dirac operators with these particular twists have properties
 similar to the Hodge--Laplacian $\D$ and the square $D^2$ of the
 untwisted Dirac operator:

 \begin{Proposition}\label{gd}
 \summary{Global Decomposition Principle}
  The restriction $D^2_{\T^{l,\,d}}\big\vert_\pi$ of the square of a
  twisted Dirac operator $D^2_{\T^{l,\,d}}$ with twisting bundle $\T^{l,\,d}
  \,:=\,(\S^lH\otimes\L^d_\circ E)(M)$ to a parallel subbundle $\pi(M)\,
  \subset\,(\Z\otimes\T^{l,\,d})(M)$ does not depend on the specific
  embedding of this subbundle and equation (\ref{tdir}) becomes in this case:
  $$
   \D_\pi
   \;=\;  D^2_{\T^{l,\,d}}\Big\vert_\pi\,+\,
    {\kappa\over 8n(n+2)}(\,l\,+\,d\,-\,n\,)\,(\,l\,-\,d\,+\,n\,+\,2\,)
  $$
 \end{Proposition}

\section{Classification of Minimal and Maximal Twists}\label{clss}

 In this section we will focus attention on the technicalities
 necessary to draw conclusions from Proposition \ref{gd}. The
 irreducible representations occurring in the twisted spinor
 representations $\Z\otimes\T^{l,\,d}$ are all of the form
 $\S^kH\otimes\M^{a,\,b}E$, where $\M^{a,\,b}E$ is the irreducible
 representation of highest weight in the tensor product $\L^a_\circ E
 \otimes\L^b_\circ E$. Alternatively we see from Weyl's construction
 of the irreducible representations of the classical matrix groups
 that $\M^{a,\,b}E$ is the common kernel of the diagonal contraction
 with the symplectic form $\s:\;\L^a_\circ E\otimes\L^b_\circ E\,
 \longrightarrow\,\L^{a-1}_\circ E\otimes\L^{b-1}_\circ E$ and
 the Pl\"ucker differential:
 $$
  \sum_\mu\,e_\mu\wedge\,\otimes\,de_\mu\inner:\quad
  \L^a_\circ E\,\otimes\,\L^b_\circ E\;\longrightarrow\;
  \L^{a+1}E\,\otimes\,\L^{b-1}_\circ E
 $$
 In particular, we will characterize the twists $\T^{l,\,d}$ with
 $\S^kH\otimes\M^{a,\,b}E\,\subset\,\Z\otimes\T^{l,\,d}$. Moreover,
 for each representation $\S^kH\otimes\M^{a,\,b}E$ in this class and
 will classify the special twists maximizing the curvature expression
 $$
  -{\kappa\over 8n(n+2)}\;(\,l\,+\,d\,-\,n\,)\,(\,l\,-\,d\,+\,n\,+\,2\,)
 $$
 of Proposition \ref{gd} for $\kappa\,>\,0$ and $\kappa\,<\,0$. This
 classification is the most important step used in the applications
 of the ideas encoded in Proposition \ref{gd}. Global questions are
 postponed to the next sections. Hence, we will deal with representations 
 of $\Sp(1)\times\Sp(n)$ only.

 \begin{Theorem}\label{admis}
 \summary{Characterization of Admissible Twists}
  A representation $\T^{l,\,d}:=\S^lH\otimes\L^d_\circ E$ with $l\geq 0$
  and $n\geq d\geq 0$ is called an admissible twist for the irreducible
  representation $\S^kH\otimes\M^{a,\,b}E$, if there exists a non--trivial,
  equivariant homomorphism from $\S^kH\otimes\M^{a,\,b}E$ to the twisted
  spinor representation $\Z\otimes\T^{l,\,d}$:
  $$
   \Hom_{\Sp(1)\times\Sp(n)}
   (\,\S^kH\otimes\M^{a,\,b}E,\;\Z\otimes \T^{l,\,d}\,)\quad\neq\quad\{0\}
  $$
  A twist $\T^{l,\,d}$ is admissible in this sense if and only if
  $\;k+a+b\,\equiv\,n+l+d\;\textrm{mod}\;2\;$ and:
  \begin{eqnarray}
   b&\leq&d\label{bottom}\\
   |\,k\,-\,l\,|\;+\;|\,a\,-\,d\,|&\leq& n\,-\,b\label{diamond}\\
   |\,n\,-\,a\,+\,b\,-\,d\,|&\leq& k\,+\,l\label{cone}
  \end{eqnarray}
  \end{Theorem}

 A simple consequence of Theorem \ref{admis} is that every irreducible
 representation $\S^kH\otimes\M^{a,\,b}E$ occurs in a twisted spinor
 representation, e.~g.~in $\Z\otimes\T^{k+n-b,\,a}$ and $\Z\otimes
 \T^{|n-a-k|,\,b}$. In fact for the twist $\T^{k+n-b,\,a}$ inequality
 (\ref{diamond}) is trivial and (\ref{cone}) needs $|n-2a+b|\,\leq\,
 |n-a|\,+\,|a-b|$. For the second twist $\T^{|n-a-k|,\,b}$ inequality
 (\ref{diamond}) follows from the distance decreasing property $||x|-|y||
 \,\leq\,|x-y|$ of the absolute value via $||-k|-|n-a-k||\,\leq\,n-a$,
 whereas (\ref{cone}) reduces to $|n-a|\,\leq\,\max\,\{n-a,\,2k-n+a\}
 \,=\,k+|n-a-k|$. These two twists are the prototype examples of
 maximal and minimal twists to be defined below.
 \pfill

 \proof
 For the proof we recall a well--known fusion rule for the
 tensor product $\L^c_\circ E\otimes\L^d_\circ E$ of the two
 irreducible $\Sp(n)$--representations $\L^c_\circ E$ and
 $\L^d_\circ E$ (cf. \cite{oni}):
 $$
  \L^c_\circ E\otimes\L^d_\circ E
  \quad=\quad
  \bigoplus_{\vbox{\hbox{$\scriptstyle\quad a+b\;\equiv\;c+d\;
  \mathrm{mod}\,2$}\vskip-2pt\hbox{$\scriptstyle \qquad a+b\;\leq\;c+d$}
  \vskip-2pt\hbox{$\scriptstyle|c-d|\;\leq\;a-b\;\leq\; 2n-c-d$}}}
  \;\;
  \M^{a,\,b}E
 $$
 Note in particular that each irreducible representation $\M^{a,\,b}E$
 occurs at most once in the tensor product $\L^c_\circ E\otimes
 \L^d_\circ E$. Using this fusion rule together with the Clebsch--Gordan
 formula for irreducible $\Sp(1)$--representations and the decomposition
 of the spinor representation $\Z$ under $\Sp(1)\times\Sp(n)$ given in
 Proposition \ref{dez} we can formally write down the decomposition
 \begin{equation}
  \bigoplus_{c=0}^n\;(\S^{n-c}H\otimes\L^c_\circ E)
   \otimes(\S^lH\otimes\L^d_\circ E)
  \;\;=\!
  \sum_{{k\geq 0\atop n\geq a\geq b\geq 0}}\!\!\!\#\;
   \mathfrak{M}_{k,\,a,\,b}(l,\,d)\cdot\S^kH\otimes\M^{a,b}E
 \end{equation}
 of $\Z\otimes\T^{l,\,d}$, where $\mathfrak{M}_{k,\,a,\,b}(l,\,d)$ is
 the set of all $n\geq c\geq 0$ satisfying the set of constraints:
 \begin{equation}\label{ineq}
  \begin{array}{ccc}
  k&\equiv&n+c+l\;\;\hbox to0pt{$\mathrm{mod}\;2$\hss}\cr
  k&\leq&n-c+l\cr k&\geq&|n-c-l|
  \end{array}
  \qquad\qquad
  \begin{array}{ccc}
  a+b&\equiv&c+d\;\;\hbox to0pt{$\mathrm{mod}\;2$\hss}\cr
  a+b&\leq&c+d\cr a-b&\geq&|c-d|\cr a-b&\leq&2n-c-d
  \end{array}
 \end{equation}

 \noindent
 It is clear from these constraints that $\mathfrak{M}_{k,\,a,\,b}
 (l,\,d)$ is empty unless $k+a+b\,\equiv\,n+l+d\;\mathrm{mod}\,2$ reflecting
 in a way the consistency of the action of $(-1,\,-1)\,\in\,\Sp(1)\times
 \Sp(n)$. In particular, $k+a+b\,\equiv\,n+l+d\;\mathrm{mod}\,2$ is a
 necessary condition for the twist $\T^{l,\,d}$ to be admissible.

 In view of this congruence we can drop one of the two constraints
 $a+b\,\equiv\,c+d\;\mathrm{mod}\,2$ or $k\,\equiv\,n+c+l\;\mathrm{mod}\,2$
 and solve the inequalities (\ref{ineq}) for $c$ to arrive after a little
 manipulation at an equivalent description of $\mathfrak{M}_{k,\,a,\,b}(l,\,d)$
 as the set of all $c\,\equiv\,a+b+d\;\mathrm{mod}\,2$ satisfying:
 \begin{equation}\label{cons}
  \max\,\{\,b+|a-d|,n-k-l\,\}\,\leq\,c\,
  \leq\,n\,-\,\max\,\{\,|k-l|,|n-a+b-d|\}
 \end{equation}

 \noindent
 Under the standing hypothesis $k+a+b\,\equiv\,n+l+d\;\mathrm{mod}\,2$
 we evidently have
 $$
  \max\,\{\,b+|a-d|,n-k-l\,\}\;\equiv\;a\,+\,b\,+\,d\;\equiv\;
  n\,-\,\max\,\{\,|k-l|,|n-a+b-d|\}\quad\mathrm{mod}\;2
 $$
 so that $\mathfrak{M}_{k,\,a,\,b}(l,\,d)$ will be non--empty if and only
 if the inequality (\ref{cons}) is consistent, because the congruence
 $c\,\equiv\,a+b+d\;\mathrm{mod}\,2$ will be fulfilled by either end of the
 resulting interval. However, the consistency condition for (\ref{cons})
 is given by four inequalities in $l,\,d$ depending of course on $k,\,a,\,b$.
 The first $n-k-l\,\leq\,n-|k-l|$ is trivial for $k,\,l\geq 0$ and the next
 two become inequalities (\ref{diamond}) and (\ref{cone}), whereas the last
 $b+|a-d|\,\leq\,n-|n-a+b-d|$ is equivalent to inequality (\ref{bottom})
 for all $b\leq a\leq n$ and $d\leq n$.
 \qed
 \pfill

 Note that if the set $\mathfrak{M}_{k,\,a,\,b}(l,\,d)$ is non-empty all
 its elements will have the same parity as $a+b+d$. Of course their number
 $\#\,\mathfrak{M}_{k,\,a,\,b}(l,\,d)$ is just the multiplicity of the
 representation $\S^kH\otimes\M^{a,\,b}E$ in $\Z\otimes\T^{l,\,d}$,
 which we will need below as index multiplicity:

 \begin{Definition}
 \summary{The Index of an Admissible Twist}
  The index of an admissible twist $\T^{l,\,d}$ for an irreducible
  representation $\S^kH\otimes\M^{a,\,b}E$ is the index multiplicity
  of $\S^kH\otimes\M^{a,\,b}E$ in the twisted spinor representation
  $\Z^{\pm}\otimes\T^{l,\,d}$:
  \begin{eqnarray*}
   \mathrm{index}\,(k,a,b;\,l,d)
   &:=&
   \dim\;\Hom_{\Sp(1)\times\Sp(n)}
   (\S^kH\otimes\M^{a,\,b}E,\,\Z^+\otimes\T^{l,\,d})\\[.5mm]
   &&\qquad-\quad
   \dim\;\Hom_{\Sp(1)\times\Sp(n)}
   (\S^kH\otimes\M^{a,\,b}E,\,\Z^-\otimes\T^{l,\,d})
  \end{eqnarray*}
  From the proof of Theorem \ref{admis} we can easily
  read off an explicit formula for this index:
  \begin{eqnarray*}
   \lefteqn{\mathrm{index}\,(k,a,b;\,l,d)\;:=\;}&&\\[.5mm]
   &&{(-1)^{a+b+d}\over 2}
    \Big(n\,+\,2\,-\,\max\,\{|k-l|,|n-a+b-d|\}
    \,-\,\max\,\{b+|a-d|,n-k-l\}\Big)
  \end{eqnarray*}
 \end{Definition}

 Although we have calculated the index multiplicity of the representation
 $\S^kH\otimes\M^{a,\,b}E$ for an arbitrary twisted spinor representation
 $\Z\otimes\T^{l,\,d}$, it will turn out below that only very few
 representations actually contribute to the index of a particular twisted
 Dirac operator. These representations are characterized by the following
 extremality condition:

 \begin{Definition}
 \summary{Minimal and Maximal Twists}
  An admissible twist $\T^{l,\,d}:=\S^lH\otimes\L^d_\circ E$ for the
  irreducible representation  $\S^kH\otimes\M^{a,b}E$ is called
  a minimal or maximal twist, if the curvature term of Proposition
  \ref{gd}, or equivalently the function $\phi(\tilde l,\tilde d)
  \;:=\;(\tilde l\,+\,\tilde d\,-\,n)\,(\tilde l\,-\,\tilde d\,+\,n\,+\,2)$,
  assumes its minimum or maximum among all admissible twists
  $\T^{\tilde l,\,\tilde d}$ in the twist $\T^{l,\,d}$.
 \end{Definition}

 To determine the index of a twisted Dirac operator in terms of the
 dimension of the eigenspaces of the operators $\D_\pi$, all we will
 further need is a classification of all minimal twists for negative
 scalar curvature $\kappa<0$ and similarly of all maximal twists for
 $\kappa>0$:

 \begin{Theorem}\label{max}
 \summary{Classification of Maximal Twists}
  All representations $\S^kH\otimes\M^{a,b}E$ with $k>0$ or $a>b$
  have unique maximal twists:
  $$
   \;\; \T^{k+n-b,\,a}\quad=\quad \S^{k+n-b}H\otimes\L^a_\circ E\qquad\qquad
   \mathrm{index}\,(k,a,b;\,k+n-b,a)\;\;=\;\;(-1)^b
  $$
  For the special representations $\M^{a,a}E$ with $k=0$ and $a=b$
  all admissible twists $\T^{n-d,d}$ with $d=a,\,\ldots,\,n$ have
  $\phi(n-d,d)=0$ and are thus automatically maximal and minimal:
  $$
   \T^{n-d,\,d}\quad\quad=\quad \S^{n-d}H\otimes\L^d_\circ E\quad\qquad\qquad
   \mathrm{index}\,(0,a,a;\,n-d,d)\quad=\quad(-1)^d
  $$
 \end{Theorem}

 \noindent
 The classification of all minimal twists splits into more cases:

 \begin{Theorem}\label{min}
 \summary{Classification of Minimal Twists}
  According to their minimal twists the irreducible representations
  $\S^kH\otimes\M^{a,b}E$ are divided into four classes. In the
  first class we have $k>(n-a)+(n-b)$ and a unique minimal twist:
  $$
   \T^{k-n+b,\,a}\quad=\quad \S^{k-n+b}H\otimes\L^a_\circ E\qquad\qquad
   \mathrm{index}\,(k,a,b;\,k-n+b,a)\;\;=\;\;(-1)^b
  $$
  In the second class with $k=(n-a)+(n-b)$ the minimal twist is
  no longer unique. All minimal twists for representations in this
  class are given by
  $$
   \T^{n-d,\,d}\quad\quad=\quad \S^{n-d}H\otimes\L^d_\circ E\quad\qquad\qquad
   \mathrm{index}\,(k,a,b;\,n-d,d)\quad=\quad(-1)^{k+d}
  $$
  with $d=b,\,\ldots,\,a$. The special representations $\M^{a,a}E$
  with $k=0$ and $a=b$ form the third class overlapping in $k=0$ and
  $a=b=n$ with the second. All admissible twists $\T^{n-d,\,d}$ with
  $d=a,\,\ldots,\,n$ for these special representations are minimal
  and maximal at the same time:
  $$
   \T^{n-d,\,d}\quad\quad=\quad \S^{n-d}H\otimes\L^d_\circ E\quad\qquad\qquad
   \mathrm{index}\,(0,a,a;\,n-d,d)\quad\;=\;\quad(-1)^d
  $$
  The remaining representations are characterized by $k<(n-a)+(n-b)$ and
  $k+(a-b)>0$. The minimal twists of the representations in this fourth
  class are all unique:
  $$
   \T^{|n-a-k|,\,b}\quad=\quad \S^{|n-a-k|}H\otimes\L^b_\circ E
   \;\;\quad\qquad\mathrm{index}\,(k,a,b;\,|n-a-k|,b)\;=\;\;(-1)^a
  $$
 \end{Theorem}

 Before proceeding to the actual proofs of Theorem \ref{max} and
 Theorem \ref{min} let us agree on some geometric terms in order to
 help intuition. The set of solutions to the inequality (\ref{diamond})
 in $(l,\,d)$--space is a ball in $L^1$--norm, i.~e.~a diamond, with
 center $(k,\,a)$ and radius $n-b$. Its right and left corner are
 thus $(k\pm(n-b),\,a)$ with $(k,\,a\pm(n-b))$ being its top and bottom
 corner. On the other hand the set of solutions to the inequality
 (\ref{cone}) is the cone $\{\,(l,\,d)\,\vert\;l+d\,\geq\,-k+n-a+b\,
 \textrm{\ and\ }\,l-d\,\geq\,-k-n+a-b\,\}$ opening diagonally to
 the right from its vertex in the point $(-k,\,n-a+b)$.

 In particular the set of solutions to both inequalities (\ref{diamond})
 and (\ref{cone}) is always a rectangle in $(l,\,d)$--space, which may
 degenerate into a straight line but always contains at least the
 points $(k+n-b,\,a)$ and $(|n-a-k|,\,b)$. Note that all corners of the
 diamond as well as the vertex of the cone and the corners of the resulting
 intersection rectangle satisfy the congruence condition $l+d\,\equiv\,
 n+k+a+b$, which consequently will care for itself below.

 Finally the level sets of the function $\phi(l,\,d)\,=\,(l+d-n)\,
 (l-d+n+2)$, which we are going to extremize, are hyperbolas with two
 diagonal axes $l+d\,=\,n$ and $l-d\,=\,-n-2$ dividing $(l,\,d)$--space
 into four quadrants. In the first quadrant with $l+d\,\geq\,n,\;
 l-d\,\geq\,-n-2$ the function $\phi\geq 0$ is positive, whereas
 it is negative in the second $l+d\,\leq\,n,\;l-d\,\geq\,-n-2$.
 Eventually we only care for points $l\geq 0$ and $n\geq d\geq 0$
 in these two quadrants.
 \pfill

 \noindent\textbf{Proof of Theorem \ref{max}:}\quad
 We already know that the right corner $(k+n-b,\,a)$ of the diamond
 always corresponds to the admissible twist $\T^{k+n-b,\,a}$ since
 $|n-2a+b|\,\leq\,|n-a|+|a-b|\,=\,n-b$. If this right corner of the
 diamond lies in the strict interior of the first quadrant, then it
 will be the unique point, where $\phi$ assumes its maximum on the
 diamond, tacitly ignoring of course third and fourth quadrant. In
 particular the twist $\T^{k+n-b,\,a}$ will be the unique maximal
 twist as soon as $k+n-b+a\,>\,n$, equivalently $k>0$ or $a>b$.

 Assuming now $k=0$ and $a=b$ we see that the top corner $(0,\,n)$ of
 the diamond coincides with the vertex of the cone. Thus the intersection
 rectangle degenerates into the face of the diamond running from its
 top corner $(0,\,n)$ to its right corner $(n-a,\,a)$. Consequently the
 admissible twists are exactly the twists $\T^{n-d,\,d}$ with $d\,=\,
 a,\,\ldots,\,n$ and $\phi(n-d,\,d)\,=\,0$. The calculation of the
 index multiplicities is left to the reader.
 \qed
 \pfill

 \noindent\textbf{Proof of Theorem \ref{min}:}\quad
 We first concentrate on the case $k\,>\,(n-a)+(n-b)$ or equivalently
 $k-n+b+a\,>\,n$, where the diamond lies completely in the strict
 interior of the first quadrant since its left corner does. With the
 axes of the level sets of $\phi$ running parallel to the faces of
 the diamond $\phi$ assumes a unique minimum on the diamond in its
 left corner. Consequently we are done once we have checked that
 $\T^{k-n+b,\,a}$ is an admissible twist. However inequality (\ref{cone})
 immediately follows from $|n-2a+b|\,\leq\,n-b\,<\,k$, which is needed
 for calculating the index multiplicity, too.

 Assuming next that $k\,=\,(n-a)+(n-b)$ the left corner of the diamond
 is the point $(n-a,\,a)$ in the first quadrant. Hence, all of the diamond
 lies in the first quadrant $\phi\geq 0$ with $\phi\,=\,0$
 only on the face from its left to its bottom corner $(2n-a-b,\,a-n+b)$.
 Note that the bottom corner fails to satisfy inequality (\ref{bottom})
 and that inequality (\ref{cone}) is satisfied by the left corner
 $(n-a,\,a)$ due to $|n-2a+b|\,\leq\,n-b\,\leq\,k$. Taking this into
 account the only admissible twists satisfying $\phi\,=\,0$ are
 exactly the twists $\T^{n-d,\,d}$ with $d\,=\,b,\,\ldots,\,a$.

 The admissible twists for the special representations $\M^{a,\,a}E$
 with $k=0$ and $a=b$ are exactly the twists $\T^{n-d,\,d}$ with $d\,=\,
 a,\,\ldots,\,n$, because the top corner $(0,\,n)$ of the diamond coincides
 with the vertex of the cone. As all these admissible twists have
 $\phi(n-d,\,d)\,=\,0$, they are all both minimal and maximal.

 Recall now that $\T^{|n-a-k|,\,b}$ is an admissible twist, because
 $|n-a|\,\leq\,k+|n-a-k|$ and $||-k|-|n-a-k||\,\leq\,n-a$ by distance
 decrease. Turning to geometry we see that the bottom corner of the
 intersection rectangle of cone and diamond will be either $(k,\,a-n+b)$
 for $k\,\geq\,n-a$ or $(n-a,\,b-k)$ for $k\,\leq\,n-a$, i.~e.~whatever
 point has larger $l$ and $d$--coordinate. In particular this bottom
 corner fails in general to satisfy inequality (\ref{bottom}) chopping
 off a triangle from the rectangle. The resulting face runs from the point
 $(|n-a-k|,\,b)$ to $(n-a+k,\,b)$ independent of whether $k\,\geq\,n-a$
 or $k\,\leq\,n-a$. Note that the geometry may become even more complicated,
 but we already know that the twist $\T^{|n-a-k|,b}$ is admissible, which
 fixes this problem as far as we need it. 

 In order to classify the minimal twists of the remaining representations
 characterized by $k\,<\,(n-a)+(n-b)$ and $k+(a-b)\,>\,0$ we observe that
 these two assumptions together are equivalent to $|n-a-k|+b\,<\,n$, so
 that the point $(|n-a-k|,\,b)$ will lie in the strict interior of the
 second quadrant. From the geometric discussion above we conclude that
 $\phi$ assumes a unique minimum in this point, because the tangents
 to the level surfaces of $\phi$ are never diagonal and horizontal
 only for $l\,=\,-1\,<\,|n-a-k|$.
 \qed

\section{Eigenvalue estimates}

 The potential applications of Proposition \ref{gd} include eigenvalue
 estimates for the Laplace and for twisted Dirac operators. The general
 procedure is described in this section and carried out in some
 particularly interesting cases. Our first example are the irreducible
 $\Sp(1)\cdot\Sp(n)$--representations $\S^rH\otimes\L^r_\circ E$ defining
 parallel subbundles in the bundle of $r$--forms (cf. \cite{sala3}).
 On these parallel subbundles we have the following lower bound for
 the spectrum of the Laplace operator.

 \begin{Proposition}
 \summary{Eigenvalue Estimate on $\S^rH\otimes\L^r_\circ E$}
  Let $(M^{4n},\,g)$ be a compact quaternionic K\"ahler manifold of
  positive scalar curvature $\kappa>0$. Then any eigenvalue $\lambda$ 
  of the Laplace operator restricted to $\S^rH\otimes\L^r_\circ E$
  satisfies
  $$
   \lambda\quad\geq\quad{r(n + 1)\over 2n(n + 2)} \, \kappa \, .
  $$
 \end{Proposition} 

 \proof
 It follows from Theorem \ref{max} that $\S^{n+r} H \otimes \L^r_\circ\,E$
 is a maximal twist for the representation  $\S^r H \otimes \L^r_\circ \,E$.
 Using Proposition \ref{gd} with $l=n+r$ and $d=r$ we obtain:
 $$
  \Delta_{\S^rH\otimes\L^r_\circ E}
  \quad=\quad
  D^2_{\T^{n+r,r}} \Big\vert_{\S^rH\otimes\L^r_\circ E}
  \; +\; {r(n + 1)\over 2n(n+2)}\,\kappa\,
  \quad\geq\quad{r(n + 1)\over 2n(n + 2)}\,\kappa \, .\qed
 $$

 \noindent
 An interesting special case is $H \otimes E = TM\otimes_\R\C$ for $r=1$,
 leading to an eigenvalue estimate for the Laplace operator on 1--forms.
 In particular, the first Betti number has to vanish. Since the differential
 of any eigenfunction of the Laplace operator is an eigenform for the same
 eigenvalue we also obtain an estimate on functions 
 (cf. \cite{alex1} and \cite{lebrun2}):

 \begin{Corollary}
 \summary{Vanishing of the First Betti Number for Positive Scalar Curvature}
  Let $(M^{4n}, \, g)$ be a compact quaternionic K\"ahler manifold of
  positive scalar curvature $\kappa>0$. Any eigenvalue $\lambda$ of 
  the Laplace operator on non--constant functions or 1--forms satisfies
  $$
   \lambda\quad\geq\quad{n+1\over 2n(n+2)} \, \kappa \, .
  $$
 \end{Corollary}
 
 \noindent  
 Replacing maximal by minimal twists to compensate the sign of the scalar
 curvature the same argument provides eigenvalue estimates on $\S^rH\otimes
 \L^r_\circ E$ on manifolds with $\kappa<0$:
 
 \begin{Proposition}\label{negative}
 \summary{Vanishing of the First Betti Number for Negative Scalar Curvature}
  Let $(M^{4n}, \, g)$ be a compact quaternionic K\"ahler manifold of
  negative scalar curvature $\kappa<0$. Then any eigenvalue $\lambda$
  of  the Laplace operator on 1--forms  satisfies:
  $$
   \lambda\quad\geq\quad{|\kappa|\over 2(n + 2)}\, .
  $$
  In particular the first Betti number has to vanish even in the case
  of negative scalar curvature.
 \end{Proposition}

 \proof
 Recall that we excluded the case $n=1$ from the very beginning in
 Definition \ref{qkh}. Since $n\geq 2$ and $r=1$ we are in the fourth
 case of Theorem \ref{min}. The unique minimal twist for $H\otimes E$
 is thus $\S^{n-2} H$ and we can apply Proposition \ref{gd} with
 $l=n-2$ and $d=0$ to obtain:
 $$
  \Delta_{H\otimes E}\;=\;D^2_{\T^{n-2,0}} \Big\vert_{H\otimes E}
  \; - \; {\kappa\over 2(n + 2)}\;\ge\;{|\kappa|\over 2(n + 2)}\, .
  \qed
 $$
 The vanishing of the first Betti number in the case of negative scalar
 curvature was also proved in \cite{horan}. In Proposition \ref{negx}
 we will prove a stronger vanishing result for the odd Betti numbers.

 As an other application we consider the Laplace operator on 2--forms
 $\L^2T^*M\otimes_\R\C$, which decompose into $\S^2H\oplus(\S^2H \otimes
 \L^2_\circ E)\oplus\S^2 E$. In the next section we will see that the
 Laplace operator may have a kernel in the sections of the parallel
 subbundle $\S^2E$. Nevertheless we have a positive lower bound on
 the other two parallel subbundles:

 \begin{Proposition}
 \summary{Eigenvalue Estimates on 2--forms}
  Let $(M^{4n}, \, g)$ be a compact quaternionic K\"ahler manifold
  of positive scalar curvature $\kappa$. Then all eigenvalues $\lambda$
  of the Laplace operator on 2--forms in $\S^2 H$ or $\S^2H\otimes
  \L^2_\circ\,E$ satisfy
  $$
   \lambda(\Delta_{\S^2H})
   \quad\geq\quad{\kappa\over 2n}\qquad\textrm{\ and\ }\qquad
   \lambda(\Delta_{\S^2H\otimes \L^2_\circ \,E})
   \quad\geq\quad{n+1\over n(n+2)}\,\kappa\ . 
  $$
 \end{Proposition} 

 The estimate for the Laplace operator on $\S^2H\subset\L^2T^*M\otimes_\R\C$ 
 was proved for the first time in \cite{alex2}. Again we have similar results 
 in the case of negative curvature. In particular, the lower bound for 
 $\Delta_{\S^2H}$ is the same as in Proposition \ref{negative}.
 \pfill

 Our next aim is to derive properties of twisted Dirac operators. For doing
 so we make the following crucial observation. If $\pi$ is any representation
 with admissible twists $\T^{l,d}$ and $\T^{{\tilde l},{\tilde d}}$ then
 we can apply Proposition \ref{gd} twice to obtain
 \begin{equation}\label{twist1}
  D^2_{\T^{l,d}}\, \Big\vert_\pi\;=\;D^2_{\T^{{\tilde l},
  {\tilde d}}}\,\Big\vert_\pi\; + \;{\kappa\over 8n(n+2)} \,\left( 
  \phi({\tilde l},\,{\tilde d})\;-\;\phi(l,\,d)\right) \ ,
 \end{equation}
 with $\phi(l,\,d)=(l+d-n)(l-d+n+2)$. We first use this observation
 to give a short proof of the eigenvalue estimate for the untwisted
 Dirac operator:

 \begin{Proposition}
 \summary{Eigenvalue Estimate for the Untwisted Dirac Operator \cite{qk1}}
  Let $(M^{4n}, \, g)$ be a compact quaternionic K\"ahler spin
  manifold of positive scalar curvature $\kappa$. Then any
  eigenvalue $\lambda$ of the untwisted Dirac operator satisfies
  $$
   \lambda^2 \quad\geq\quad {n+3\over n+2} \, {\kappa\over 4} \ .
  $$
 \end{Proposition}

 \proof
 According to Proposition \ref{dez} the spinor bundle decomposes
 into the parallel subbundles $\Z\,=\,\oplus^n_{r=0}\;\Z_r$ with
 $\Z_r\,=\,\S^rH\otimes\L^{n-r}_\circ E$. To estimate the square
 of the Dirac operator on $\S^rH\otimes\L^{n-r}_\circ E$ we observe
 that the unique maximal twist for $\S^rH\otimes\L^{n-r}_\circ E$
 is $\T^{n+r, n-r}$ and for $l=d=0$ and $\tilde l=n+r,\,\tilde d=n-r$
 equation (\ref{twist1}) reads:
 $$
  D^2 \Big\vert_{\Z_r}\quad=\quad D^2_{\T^{n+r, n-r}}\Big\vert_{\Z_r}
  \; + \;{\kappa\over 8n(n+2)} \,\Big(\,n(2r+n+2)\,+\,n(n+2)\,\Big)
  \quad\geq\quad{n+2+r\over n+2}\,{\kappa\over 4}\;.
 $$
 Consequently some hypothetical eigenspinor $\phi\in\Gamma(\Z)$ of 
 $D^2$ with eigenvalue $\lambda^2\;<\;{n+3\over n+2}{\kappa\over4}$ 
 would have to be localized in the subbundle $\Z_0\subset\Z$. But
 the Dirac operator on a manifold of positive scalar curvature has
 trivial kernel so that $D\phi\in\Gamma(\Z_1)$ would be a nontrivial
 eigenspinor for $D^2$ again with eigenvalue $\lambda^2$ in
 contradiction to the estimate for $\Z_1$.
 \qed
 \pfill
  
 We now use equation (\ref{twist1}) for describing the kernels of twisted
 Dirac operators in the case of positive scalar curvature. If $\pi$ is 
 any representations which contributes to the kernel of $D^2_{\T^{l,d}}$ 
 then $\T^{l,d}$ has to be a maximal twist for $\pi$. In fact equation
 (\ref{twist1}) implies that $D^2_{\T^{l,\,d}}$ is positive on $\pi$ 
 as soon as there is another admissible twist $\T^{\tilde l,
 \,\tilde d}$ for $\pi$ with $\phi({\tilde l}, {\tilde d})>\phi(l, d)$.
 From this remark and Proposition \ref{gd} we conclude in the case of
 positive scalar curvature
 \begin{equation}\label{kernel}
  \ker(D^2_{\T^{l,d}}) \; = \; \bigoplus_\pi 
  \, \ker\left(\Delta_\pi \, - \, {\kappa\over 8n(n+2)} \phi(l, d)\right) 
 \end{equation}
 where the sum is over all representations $\pi$ for which $\T^{l,\,d}$
 is a maximal twist. Since ${\kappa\over 8n(n+2)}\,\phi(l, d)$ is the
 smallest possible eigenvalue of the operator $\Delta_\pi$ equation
 (\ref{kernel}) is in essence a decomposition of $\ker(D^2_{\T^{l,d}})$
 into a sum of minimal eigenspaces for the operators $\Delta_\pi$.

 If $\T^{l,\,d}$ is a maximal twist for a representation $\pi$ then Theorem
 \ref{max} also provides us with the information whether $\pi$ occurs in
 $\Z^+ \otimes \T^{l,\,d} $ or in $\Z^- \otimes \T^{l,\,d}$. Hence a
 corollary of equation~(\ref{kernel}) is a formula for the index of the
 twisted Dirac operator $D_{\T^{l,\,d}}$ in terms of dimensions of certain
 minimal eigenspaces. We will describe this in two examples:

 \begin{Proposition}
 \nosummary
  Let $(M^{4n}, \, g)$ be a compact quaternionic K\"ahler
  manifold of positive scalar curvature $\kappa>0$, then:
  $$
   \ker\left(D^2_{\T^{l,d}}\right) \; = \; \{0\}
   \qquad \textrm{\ for\ } \qquad l \, + \, d \, < \,  n \ .
  $$
 \end{Proposition}

 \proof
 All maximal twists $\T^{l,d}$ satisfy $l+d\geq n$ by Theorem \ref{max}.
 \qed

 \pfill\noindent
 An immediate consequence of this proposition is the vanishing
 of the index $\ind(D_{\T^{l,\,d}})$ for $l+d<n$. This was also
 proved in \cite{lebrun} by using the Akizuki--Nakano vanishing
 theorem on the twistor space.
 For the second example we consider the twisted Dirac operator
 $D_{\T^{n+2,\,0}}$. It easily follows from Theorem \ref{max} that
 $\S^2H$ is the only representation with maximal twist $\T^{n+2,\,0}$:

 \begin{Proposition}
 \summary{Killing Vector Fields}
  On every compact quaternionic K\"ahler manifold $(M^{4n}, g)$
  of positive scalar curvature $\kappa$ we have:
  $$
   \ker\left(D^2_{\T^{n+2,0}}\right)\; = \; 
   \ker\left(\Delta_{\S^2 H} \, - \, {\kappa\over 2n}\right) \ .
  $$
 \end{Proposition}

 \noindent
 The index of $D_{\T^{n+2,0}}$ equals the dimension of the isometry group
 of $(M, g)$ (cf. \cite{sala}). But since $\S^2 H$ is the only representation
 contributing to $\ker(D^2_{\T^{n+2,0}})$ the index is just the dimension of
 the minimal eigenspace of $\Delta_{\S^2 H}$. In fact, there is an explicit
 isomorphism from the space of Killing vector fields to $\S^2 H$
 (cf.~\cite{alex2}). It is given by projecting the covariant derivative
 of a Killing vector field onto its component in $\S^2H\subset\L^2T^*M
 \otimes_\R\C$.

\section{Harmonic forms and Betti numbers}

 This section contains the most important application of Proposition
 \ref{gd}. We will determine which parallel subbundles of the differential
 forms may carry harmonic forms and thus prove vanishing theorems for
 Betti numbers both for positive and negative scalar curvature.
 These results will lead to quaternionic K\"ahler analogues of
 the weak and strong Lefschetz theorem in K\"ahler geometry.
 Recall that the weak Lefschetz theorem for K\"ahler manifolds $M$
 states the inequality $b_k\le b_{k+2}$ of the Betti numbers for
 $k<{1\over2}\dim\,M$, whereas the strong Lefschetz theorem asserts
 that the wedge product with the parallel 2--form descends to an
 injective map of the cohomology $H^k(M,\,\R)\,\longrightarrow
 \,H^{k+2}(M,\,\R)$.

 \begin{Proposition}\label{harmonic}
 \summary{Representations and Harmonic Forms}
  Let $(M^{4n}, \, g)$ be a compact quaternionic K\"ahler manifold of
  scalar curvature $\kappa\neq 0$ and let $\pi$ be an irreducible
  representation of $\Sp(1)\cdot\Sp(n)$ occurring in the forms
  $\L^\bullet(H\otimes E)$:
  $$
   \Hom_{\Sp(1)\cdot\Sp(n)}(\pi, \L^\bullet(H\otimes E))\;\neq\;\{0\}
  $$
  If the scalar curvature is positive then $\ker(\Delta_\pi)
  \,=\,\{0\}$ unless $\pi\,=\,\M^{a,\,a}E$ for some $a$ with $n\geq a\geq 0$.
  Similarly if the scalar curvature is negative then
  $\ker(\Delta_\pi)\,=\,\{0\}$ unless either $\pi\,=\,\M^{a,\,a}E$
  as before or $\pi$ is a representation of the form $\pi\,=\,\S^{2n-a-b}
  H\otimes\M^{a,\,b}E$ with
  $n\geq a\geq b\geq 0$.
 \end{Proposition}

 Although the representations $\S^{2n-a-b}H\otimes\M^{a,\,b}E$ form
 a larger class of representations they are still rather special
 among all the representations occurring in the forms. The appearance
 of these exceptional representations potentially carrying harmonic
 forms could have been foreseen from the difficulties encountered
 in the attempt to push Kraines original strong Lefschetz theorem
 (\cite{krain}) for quaternionic K\"ahler manifolds beyond degree
 $n$. In higher degrees the given proofs fail precisely for these
 representations. It follows from Proposition \ref{harmonic} that
 this problem is absent in the positive scalar curvature case. 
 \pfill
 \proof
 For any manifold of even dimension the bundle of exterior forms is
 the tensor product of the spinor bundle with itself. The decomposition
 of $\Z$ given in Proposition \ref{dez} implies:
 $$
  \L^\bullet(H\otimes E) \quad=\quad \Z\otimes\Z
  \quad=\quad\bigoplus_{r=0}^n\;\Z\otimes\T^{r,\,n-r}\,.
 $$
 In particular, a representation $\pi$ occurs in the forms if and only
 if it occurs in a twisted spinor bundle $\Z\otimes\T^{r,\,n-r}$ for some $r$
 with $n\geq r\geq 0$. It is consequently of the form $\pi=\S^kH\otimes
 \M^{a,\,b}E$ for suitable $k\geq 0$ and $n\geq a\geq b\geq 0$. In this
 situation Proposition \ref{gd} becomes:
 $$
  \Delta\Big\vert_{\pi}
  \quad=\quad\Delta_\pi
  \quad=\quad D^2_{\T^{r,\,n-r}}\Big\vert_{\pi}
 $$
 A harmonic form in the parallel subbundle determined by $\pi$ is thus
 identified with an harmonic twisted spinor for the twist $\T^{r,\,n-r}$.
 However, we have already expressed the kernel of the twisted Dirac
 operators $D^2_{\T^{r,\,n-r}}$ in formula (\ref{kernel}) at least
 for positive scalar curvature.

 The point in this formula is of course that only those representations
 $\pi$ may contribute to the kernel of the twisted Dirac operator
 $D^2_{\T^{r,\,n-r}}$, for which the twist $\T^{r,\,n-r}$ is a maximal
 twist. Replacing maximal by minimal twists the same argument applies in
 the case of negative scalar curvature and we conclude that a representation
 $\pi$ may carry harmonic forms in the case of negative or positive scalar
 curvature if and only if it has a minimal or maximal twist respectively
 of the form $\T^{r,\,n-r}$ for some $r$ with $n\geq r\geq 0$. A look at the
 classification of maximal and minimal twists in Theorems \ref{max}
 and \ref{min} completes the proof.
 \qed
 \pfill
 We now want to point out a remarkable property of minimal and
 maximal twists: If a twist $\T^{l,\,d}$ is minimal or maximal for a 
 representation $\pi$ then $\pi$ always occurs with multiplicity one in the 
 twisted spinor representation $\Z\otimes\T^{l,\,d}$. Although this property 
 seems very natural it is obtained only as a corollary of the 
 calculation of the index multiplicities in Theorems \ref{max} and \ref{min} 
 using all the rather technical calculations of that section. Surely it is 
 tempting to search for a direct argument providing better insight into the 
 nature of this property.

 For us this property is very convenient counting the total multiplicity
 of those representations $\pi$ in the differential forms, which may carry
 harmonic forms. In fact for any representation $\pi$ this total
 multiplicity is given by:
 \begin{equation}\label{count}
  \dim\,\Hom_{\Sp(1)\cdot\Sp(n)}(\,\pi,\,\L^\bullet(H\otimes E)\,)
  \quad=\quad\sum_{r=0}^n\;
  \dim\,\Hom_{\Sp(1)\cdot\Sp(n)}(\,\pi,\,\Z\otimes\T^{r,\,n-r}\,)\;.
 \end{equation}
 However, in the course of the proof of Theorem \ref{harmonic} we
 characterized the representations $\pi$ potentially carrying harmonic
 forms in negative or positive scalar curvature by their property of
 having a minimal or maximal twist respectively of the form
 $\T^{r,\,n-r},\,n\geq r\geq 0$. For such a representation $\pi$
 a twist of the form $\T^{\tilde r,\,n-\tilde r}$ is minimal or maximal
 respectively if and only if it is admissible, because in this case
 $\phi(r,n-r)\,=\,0\,=\,\phi(\tilde r,n-\tilde r)$.

 Consequently for any representation $\pi$ which may carry harmonic forms
 the summands on the right hand side of equation (\ref{count}) are all either
 0 or 1 and the total multiplicity of $\pi$ in the differential forms is just
 the number of different minimal or maximal twists respectively. This number
 is easily read off from Theorems \ref{max} and \ref{min} and is part of the
 following lemma:

 \begin{Lemma}\label{lef}
 \summary{Embeddings of Harmonic Forms}
  The representation $\pi=\M^{a,\,a}E,\,n\geq a\geq 0$, occurs $n-a+1$
  times in the forms: it occurs with multiplicity one in the forms
  of degree $2a,\,2a+4,\,2a+8,\,\ldots,\,4n-2a$. Similarly the
  representation $\pi=\S^{2n-a-b}H\otimes\M^{a,\,b}E,\,n\geq a\geq b\geq 0$,
  occurs in the forms of degree $2n-a+b,\,2n-a+b+2,\,2n-a+b+4,\,\ldots,\,
  2n+a-b$ with multiplicity 1 and $a-b+1$ times in total.
 \end{Lemma}

 \pfill
 \proof
 We have already calculated the total multiplicity of the representations
 $\M^{a,\,a}E$ and $\S^{2n-a-b}H\otimes\M^{a,\,b}E$ in the differential
 forms so that it is sufficient to prove the existence of embeddings of
 these representations into the forms of the claimed degrees. First let
 us recall the well known general decomposition of the exterior forms
 $\L^k(\,H\otimes E\,)$ into Schur functors
 $$
  \L^k(\,H\otimes E\,)
  \quad=\quad\bigoplus_\Y \W_\Y H\otimes\W_{\overline{\Y}}E
 $$
 where the sum is over all Young tableaus $\Y$ of size $|\Y|=k$ and
 $\overline{\Y}$ denotes the conjugated Young tableau (\cite{fh}). All
 Schur functors have two preferred realizations as the images
 of Schur symmetrizers in iterated tensor products. Specifying the Young
 tableau $\Y$ either by the length of its rows $(r_1,\,r_2,\,\ldots,r_{c_1})$
 or of its columns $(c_1,\,c_2,\,\ldots,\,c_{r_1})$ satisfying
 $r_1\geq r_2\geq \ldots\geq r_{c_1}$ and $c_1\geq c_2\geq\ldots\geq
 c_{r_1}$ these two preferred realizations of the Schur functors
 \begin{eqnarray*}
  \W_\Y H &\subset&
   \;\;\L^{c_1}H\;\;\otimes\;\;\L^{c_2}H\;\;\otimes\;\;
   \ldots\;\;\otimes\;\;\L^{c_{r_1}}H\\
  \W_\Y E &\subset&
   \S^{r_1}E \otimes \S^{r_2}E\otimes\;\;\ldots\;\;\otimes\S^{r_{c_1}}E
 \end{eqnarray*}
 are given as the intersection of the kernels of all possible Pl\"ucker
 differentials. In our case all Schur functors in $H$
 corresponding to Young tableaus of more than two rows vanish and since
 $\L^2H\cong\C$ is trivial the Schur functor in $H$ for the Young tableau
 of size $k$ with two rows $(k-s,\,s)$ is equivalent to $\S^{k-2s}H$:
 $$
  \L^k(H\otimes E)\quad=\quad
  \bigoplus_{s=0}^{\lfloor{k\over2}\rfloor}\;\;\S^{k-2s}H
  \otimes\W_{\overline{(k-s,\,s)}}E\,.
 $$
 Conjugation of Young tableaus is defined by exchanging rows and
 columns. Conjugated to the Young tableau with two rows $(k-s,\,s)$
 is the tableau with two columns $\overline{(k-s,\,s)}$. Thus
 $\W_{\overline{(k-s,\,s)}}E$ can be defined as the kernel
 of the Pl\"ucker differential:
 $$
  \sum_\mu e_\mu\wedge\otimes de_\mu\inner:
  \quad\L^{k-s}E\otimes\L^s E\;\longrightarrow\;
  \L^{k-s+1}E\otimes\L^{s-1}E\,.
 $$
 From Weyl's construction of the representation $\M^{a,\,b}E$ as the
 intersection of the kernel of the Pl\"ucker differential $\L^a_\circ E
 \otimes\L^b_\circ E\longrightarrow\L^{a+1}E\otimes\L^{b-1}_\circ E$ with
 the kernel of the diagonal contraction with the symplectic form we see
 that $\M^{a,\,a}E\subset\W_{\overline{(a,\,a)}}E$. Consider now the map
 $$
  \Omega:\quad\L^aE\otimes\L^bE\;\longrightarrow\;\L^{a+2}E\otimes\L^{b+2}E
 $$
 defined by
 $$
  \Omega:=\sum_{\mu,\nu}\left(\,
  de^\b_\mu\wedge de^\b_\nu\wedge\,\otimes\,e_\mu\wedge e_\nu\wedge\;+\;
  de^\b_\mu\wedge e_\mu\wedge\,\otimes\,de^\b_\nu\wedge e_\nu\wedge\,\right)\,,
 $$
 which curiously enough commutes with the Pl\"ucker differential.
 Consequently we may extend the above embedding to a chain of
 $\Sp(n)$--equivariant linear maps:
 $$
  \M^{a,\,a}E\;\longrightarrow\;\W_{\overline{(a,\,a)}}E\;
  \stackrel\Omega\longrightarrow\;\W_{\overline{(a+2,\,a+2)}}E\;
  \stackrel\Omega\longrightarrow\;\ldots\;\stackrel\Omega\longrightarrow\;
  \W_{\overline{(2n-a,\,2n-a)}}E\,.
 $$
 Explicit calculation shows that $\Omega^{n-a}=(2n-2a+1)!\,(\star\otimes
 \star)$ on $\M^{a,\,a}E$, where $\star$ denotes the Hodge isomorphism $\L^aE
 \longrightarrow\L^{2n-a}E$. Hence $\M^{a,\,a}E$ embeds into all the Schur
 functors $\W_{\overline{(a+2s,\,a+2s)}}E$ with $n-a\geq s\geq 0$ and further
 into the forms $\L^{2a+4s}(\,H\otimes E\,)$ of degree $2a+4s$ with $n-a
 \geq s\geq 0$. The appearance of the map $\Omega$ is by no means an accident,
 it can be shown that it corresponds exactly to the wedge product with the
 parallel Kraines form $\Omega$ on the level of forms.

 The construction of the different embeddings of the representations
 $\S^{2n-a-b}H\otimes\M^{a,\,b}E$ is simpler, although it is a dead end
 to start with the inclusion $\M^{a,\,b}E\subset\W_{\overline{(a,\,b)}}E$.
 Instead we have to use the Hodge isomorphism $(\star\otimes 1):\;\L^aE
 \otimes\L^bE\longrightarrow\L^{2n-a}E\otimes\L^bE$, which interchanges
 in a sense the roles of the Pl\"ucker differential and the diagonal
 contraction with the symplectic form. The Hodge isomorphism can be
 extended to a chain of maps
 $$
  \M^{a,\,b}E\;\longrightarrow\;\L^{2n-a}E\otimes \L^bE
  \;\stackrel\sigma\longrightarrow\;\L^{2n-a+1}E\otimes\L^{b+1}E
  \;\stackrel\sigma\longrightarrow\;\ldots
  \;\stackrel\sigma\longrightarrow\;\L^{2n-b}E\otimes\L^aE
 $$
 using the diagonal multiplication $\sigma$ with the symplectic form.
 Since diagonal contraction and multiplication with the symplectic form
 generate an $\mathfrak{sl}_2$--algebra of operators the final map
 $\M^{a,\,b}E\longrightarrow\L^{2n-b}E\otimes\L^aE$ is injective
 and maps into the kernel of $\sigma$. In addition the commutator
 relations between the Pl\"ucker differential and $\sigma$ imply that
 $\M^{a,\,b}E$ is mapped into the kernel $\W_{\overline{(2n-a+s,\,b+s)}}E$
 of the Pl\"ucker differential at each step, so that
 $$
  \S^{2n-a-b}H\otimes\M^{a,\,b}E\;\longrightarrow\;
  \S^{2n-a-b}H\otimes\W_{\overline{(2n-a+s,\,b+s)}}E\;
  \stackrel\subset\longrightarrow\;\L^{2n-a+b+2s}(\,H\otimes E\,)
 $$
 embeds into the forms of degree $2n-a+b+2s$ for all $a-b\geq s\geq 0$.
 \qed

 \begin{Remark}
 \summary{Strong Lefschetz Theorems}
  In the course of the proof of Lemma \ref{lef} we have sketched a
  proof of the strong Lefschetz Theorem for quaternionic K\"ahler
  manifolds of positive scalar curvature. The wedge product with the
  parallel Kraines form $\Omega$ is injective on the forms of type
  $\M^{a,\,a}E$ in all degrees $k<{1\over2}\dim\,M$ and hence descends
  to an injective map of the cohomology $H^{k}(M,\,\R)\,\longrightarrow
  \,H^{k+4}(M,\,\R)$.

  A completely different argument can be given to show that the wedge
  product with the Kraines form is injective on forms of type
  $\S^{2n-a-b}H\otimes\M^{a,\,b}E$ in degrees $k<{1\over2}
  \dim\,M\,-\,1$, too. In contrast to the positive scalar curvature case
  however, the decomposition of the cohomology given in Proposition
  \ref{harmonic} for quaternionic manifolds of negative scalar curvature
  is finer than the decomposition into primitive cohomologies with
  respect to the Kraines form.
 \end{Remark}

 \noindent
 The weak Lefschetz Theorem for quaternionic K\"ahler manifolds of
 positive scalar curvature was proved by S.~Salamon (cf.~\cite{sala})
 by analyzing the cohomology of the twistor space. Applying Proposition
 \ref{harmonic} in combination with Lemma \ref{lef} we get a more
 explicit version of this result:

 \begin{Proposition}\label{positive}
 \summary{Weak Lefschetz Theorem for Positive Scalar Curvature}
  Let $(M^{4n},\,g)$ be a compact quaternionic K\"ahler
  manifold of positive scalar curvature $\kappa>0$. Its
  Betti numbers $b_k$ satisfy for all $0\leq k\leq n$ the
  following relations:
  $$\begin{array}{lll}
    \quad b_{2k+1} & = & 0\,, \\[1.5ex]
    \quad\; b_{2k} & = & \sum_{\nu=0}^{\lfloor{k\over2}
     \rfloor}\,\dim\,(\ker\Delta_{\M^{k-2\nu,\,k-2\nu}E})\,,\\[1.5ex]
    b_{2k}\;-\;b_{2k-4} & = & \dim\,(\ker\Delta_{\M^{k,\,k}E})\quad\ge\quad0\,.
  \end{array}$$
 \end{Proposition}

 \proof
 For a compact quaternionic K\"ahler manifold of positive scalar curvature
 it follows from Proposition~\ref{harmonic} that the only representations 
 potentially carrying harmonic forms are $\M^{a,\,a}E$ with $n\geq a\geq 0$.
 But according to Lemma \ref{lef} all these representations embed into
 forms of even degree, i.~e.~all odd Betti numbers necessarily vanish.
 Moreover the representations $\M^{a,\,a}E$ occur in the forms of degree
 $2k$ if and only if $a=k,\,k-2,\,\ldots$ and in this case they occur
 with multiplicity one.
 \qed
 \pfill

 \begin{Remark}
 \summary{Associated Twistor Space and 3--Sasakian Manifold \cite{galsal}}
  Let $\mathcal S$ be the 3--Sasakian manifold and $\mathcal Z$ the
  twistor space associated with the quaternionic K\"ahler manifold $M^{4n}$.
  The dimension of $\ker\Delta_{\M^{k,\,k}E}$ can be reinterpreted as
  the dimension of the cohomology of $\mathcal{S}$ and as the dimension
  of the primitive cohomology group of $\mathcal{Z}$:
  $$
   \dim(\ker\Delta_{\M^{k,\,k}E})\;=\; b_{2k}(\mathcal S)
   \;=\; b_{2k}(\mathcal Z)\;-\; b_{2k-2}(\mathcal Z)
   \qquad\qquad k\leq n\;\ .
  $$
 \end{Remark}

 As an immediate consequence of Proposition \ref{positive} we
 obtain a result of S.~Salamon and C.~LeBrun (cf.~\cite{lebrun})
 on the index of the twisted Dirac operator $D_{\T^{l,\,d}}$
 with $l+d\,=\,n$:

 \begin{Corollary}\label{cor}
 \summary{Index of Twisted Dirac Operators and Betti Numbers}
  Let $(M^{4n}, \, g)$ be a compact quaternionic K\"ahler
  manifold of positive scalar curvature $\kappa>0$:
  $$\begin{array}{lll}
   \;\;\ker (D^2_{\T^{n-d,\,d}}) & = & \bigoplus_{a\leq d} \,
    \ker(\Delta_{\M^{a,\,a}E})\,,\\[1.5ex]
   \dim \ker (D^2_{\T^{n-d,\,d}}) & = &
    \quad b_{2d} \; + \; b_{2d-2} \,,\\[1.5ex]
   \;\;\ind(D_{\T^{n-d,\,d}}) & = & (-1)^d\big(b_{2d}\;+\;b_{2d-2}\big)\,. 
  \end{array}$$
 \end{Corollary}

 \proof
 We already observed in formula (\ref{kernel}) that in the case of positive
 scalar curvature a representation $\pi$ may contribute to the kernel of a
 twisted Dirac operator $D^2_{\T^{l,\,d}}$ only if the twist $\T^{l,\,d}$
 is maximal for $\pi$. On the other hand the twisted spinor representation
 $\Z\otimes\T^{n-d,\,d}$ occurs in the forms so that a representation $\pi$
 contributes to the kernel of $D^2_{\T^{n-d,\,d}}$ if and only if it carries
 harmonic forms, i.~e.~$\pi$ must be one of the representations $\M^{a,\,a}E$
 for some $a$ with $n\geq a\geq 0$. From equation (\ref{bottom}) of Theorem
 \ref{admis} it is evident that $\pi=\M^{a,\,a}E$ occurs in $\Z\otimes
 \T^{n-d,\,d}$ if and only if $a\leq d$. Consequently Proposition
 \ref{positive} provides the expression for the dimension of the
 kernel of $D_{\T^{n-d,\,d}}$ in terms of Betti numbers.
 \qed
 \pfill
 In dealing with quaternionic K\"ahler manifolds of negative scalar
 curvature it is convenient to decompose their cohomology into two
 direct summands with quite different behavior:

 \begin{Definition}
 \summary{\ $\sp(1)$--Invariant and Exceptional Cohomology}
  Let $(M^{4n},\,g)$ be a compact quaternionic K\"ahler manifold of
  negative scalar curvature. According to Proposition \ref{harmonic}
  two different series of representations contribute to harmonic forms
  on $M$, namely $\M^{a,\,a}E,\,n\geq a\geq 0$ and $\S^{2n-a-b}H\otimes
  \M^{a,\,b}E,\,n\geq a\geq b\geq 0$. In particular the de Rham cohomology
  of $M$ splits into the direct sum
  $$
   H^\bullet_{dR}(\,M,\,\C\,)\quad=\quad H^\bullet_{\sp(1)}(\,M,\,\C\,)\;
   \oplus\;H^\bullet_\x(\,M,\,\C\,)
  $$
  of its $\sp(1)$--invariant cohomology $H^\bullet_{\sp(1)}(\,M,\,\C\,)$,
  which is the sum of all isotypical components corresponding to the
  representations $\M^{a,\,a}E,\,n\geq a\geq 0$, and its exceptional
  cohomology $H^\bullet_\x(\,M,\,\C\,)$, which is the direct sum of
  all isotypical components corresponding to the remaining representations
  $\S^{2n-a-b}H\otimes\M^{a,\,b}E,\,n\geq a \geq b\geq 0,\,b\neq n$.
 \end{Definition} 

 Because the curvature tensor of $M$ is $\sp(1)$--invariant the same is
 true for all its characteristic classes. Moreover $H^\bullet_{\sp(1)}
 (\,M,\,\C\,)$ is closed under multiplication and the decomposition of the
 de~Rham--cohomology into $\sp(1)$--invariant and exceptional cohomology
 is respected by the induced modul structure. A deeper analysis of the
 ring structure of the cohomology ring of $M$ will be given in a
 forthcoming paper (cf.~\cite{gregor}).

 As a final application of the ideas developed in this article we combine
 Proposition~\ref{harmonic} and Lemma~\ref{lef} to obtain new information
 on the Betti numbers of compact quaternionic K\"ahler manifolds of
 negative scalar curvature.
 
 \begin{Proposition}\label{negx}
 \summary{Weak Lefschetz Theorem for Negative Scalar Curvature}
  Let $(M^{4n}, \, g)$ be a compact quaternionic K\"ahler manifold of
  negative scalar curvature $\kappa<0$. Its $\sp(1)$--invariant and
  exceptional Betti numbers $b_{\sp(1),\,k}$ and $b_{\x,\,k}$ satisfy:
  $$\begin{array}{llll}
   b_{\sp(1),\,k} & =    & 0 \qquad\qquad &
    \textrm{\ for }\; k\textrm{\ odd}\,,\\[.5ex]
   b_{\x,\,k}     & =    & 0              &
    \textrm{\ for }\; k\;\le\;n-1\,,\\[.5ex]
   b_{\sp(1),\,k} & \leq & b_{\sp(1),\,k+4} &
    \textrm{\ for }\; k\;\leq\;2n-2\,,\\[.5ex]
   b_{\x,\,k}     & \leq & b_{\x,\,k+2}     &
    \textrm{\ for }\; k\;\leq\;2n-1\,.
  \end{array}$$
  In particular, its Betti numbers $b_k=b_{\sp(1),\,k}+b_{\x,\,k}$ satisfy:
  $$\begin{array}{llll}
   b_{2k+1} & = \;   & 0\qquad\qquad &
    \textrm{\ for  }\; 2k+1\;\le\; n-1\,,\\[.5ex]
   b_{k} & \leq\;  & b_{k+2} &
    \textrm{\ for odd }\,k\;\leq\; 2n-1\,,\\[.5ex]
   b_{k} & \leq\; & b_{k+4}  &
    \textrm{\ for  }\qquad k\;\leq\;2n-2\,.
  \end{array}$$
   \end{Proposition}
 
 \proof
 Since the $\sp(1)$--invariant Betti numbers correspond by definition
 to the representations  $\M^{a,\,a}E,\,n\geq a\geq 0$, they have the
 same properties as Betti numbers of a quaternionic K\"ahler manifolds
 of positve scalar curvature given in Proposition \ref{positive}.

 It follows from Lemma \ref{lef} that the remaining representations
 $\S^{2n-a-b}H\otimes\M^{a,\,b}E$ with $\,n\geq a\geq b\geq 0$ and
 $b\neq n$ corresponding to the exceptional Betti numbers embed into
 forms of degree $2n-a+b,\,2n-a+b+2,\,\ldots,\,2n+a-b$. For $a\not
 \equiv b\;\textrm{mod}\,2$ these embeddings give rise to harmonic
 forms of odd degree. Nevertheless the odd Betti numbers of degree
 less than $n$ have to vanish because of $2n-a+b\geq n$.
 \qed
\end{document}